\documentclass[11pt]{article}
\usepackage{amsthm, amsmath , amsfonts,  eucal } 
\usepackage[all]{xy}

\usepackage{graphics}
\makeatletter
\renewcommand{\@oddhead}{\it 
  Fass\`o \& Sansonetto:
  Symplectization of fully-Hamiltonian almost-symplectic systems 
\hfill \thepage}
\makeatother

\theoremstyle{plain}
\newtheorem{theorem}{\bf Theorem}
\newtheorem{proposition}[theorem]{\bf Proposition}
\newtheorem{lemma}[theorem]{\bf Lemma}
\newtheorem{definition}[theorem]{\bf Definition}
\newtheorem{corollary}[theorem]{\bf Corollary}
\newtheorem{remark}[theorem]{\bf Remark}
\newtheorem{remarks}[theorem]{\bf Remarks}

\usepackage{titlesec}
\titleformat{\section}
  {\normalfont\large\bfseries}{\thesection}{0.5em}{}
\newcommand{\addperiod}[1]{#1.}
\titleformat{\subsection}[runin]
  {\normalfont\bfseries}{\thesubsection}{0.25em}{\addperiod}

\begin{document}


\renewcommand\em{\it}   

\newcommand\imm{\iota}

\newcommand\fcomment[1]{\marginpar{\tiny #1}}
\newcommand\unity{\mathbb I}

\newcommand\bList{
\begin{list}{}{\leftmargin2em\labelwidth1.2em\labelsep.5em\itemindent0em
\topsep0ex\itemsep-.8ex} }
\newcommand\eList{\end{list}}


\newcommand{\beq}[1]{\begin{equation}\label{#1}}
\newcommand{\eeq}{\end{equation}}
\newcommand{\for}[1]{(\ref{#1})}

\renewcommand{\a}{\alpha}
\renewcommand{\P}{\Phi}
\newcommand{\s}{\sigma}
\renewcommand{\o}{\omega}

\newcommand{\cA}{\mathcal{A}}                 
\newcommand{\cB}{\mathcal{B}}
\newcommand{\cC}{\mathcal{C}}                 
\newcommand{\cD}{\mathcal{D}}                 
\newcommand{\cF}{\mathcal{F}}                 
\newcommand{\cG}{\mathcal{G}}                 
\newcommand{\cH}{\mathcal{H}}                 
\newcommand{\cL}{\mathcal{L}}                 
\newcommand{\cO}{\mathcal{O}}                 
\newcommand{\cP}{\mathcal{P}}                 

\newcommand{\Sp}[1]{ \mathrm{Sp}(#1)}

\newcommand{\toro}{\mathbb{T}}
\newcommand{\reali}{\mathbb{R}}
\newcommand{\razionali}{\mathbb{Q}}
\newcommand{\interi}{\mathbb{Z}}
\newcommand{\tenne}{\toro^n}
\newcommand{\renne}{\reali^n}

\newcommand{\der}[2]{\frac{\partial#1}{\partial#2}}
\newcommand{\const}{\mathrm{const}} 
\newcommand\rank{\mathrm{rank\,}}

\newcommand\saa{\s^\mathrm{aa}}

\title{\bf \parskip=0pt
Symplectization of certain Hamiltonian systems\\
in fibered almost-symplectic manifolds 
}

\author{
Francesco Fass\`o\footnote{
Universit\`a di Padova,
Dipartimento di Matematica ``Tullio Levi-Civita'',
Via Trieste 63, 35121 Padova, Italy.
{\tt (E-mail: fasso@math.unipd.it)}
}
\ and 
Nicola Sansonetto\footnote{
Universit\`a di Verona, 
Dipartimento di Informatica, 
Strada le Grazie 15, 37134 Verona, Italy.
 }
}

\vskip 1truecm
\date{\small (\today)}

\maketitle

\baselineskip=14pt   

{\small
\begin{abstract}
\noindent
There is an important difference between Hamiltonian-like vector
fields in an almost-symplectic manifold $(M,\s)$, compared to the
standard case of a symplectic manifold: in the almost-symplectic
case, a vector field such that the contraction $i_X\s$ is closed
need not be a symmetry of $\s$. We thus call {\it partially-Hamiltonian}
those vector fields which have the former property and
{\it fully-Hamiltonian} those which have both properties.
We consider $2n$-dimensional almost symplectic manifolds with a fibration
$\pi:M\to B$ by Lagrangian tori. Trivially, all vertical
partially-Hamiltonian vector fields are fully-Hamiltonian. We
investigate the existence and the properties of non-vertical
fully-Hamiltonian vector fields. We show that this class is non-empty, 
but under certain genericity conditions that involve the Fourier
spectrum of their Hamiltonian, these vector fields reduce, under
a (possibly only semi-globally defined, if $n\ge4$)
torus action, to families of
standard symplectic-Hamiltonian vector fields.

\noindent
\end{abstract}
}

{\small
}

\section{Introduction}

\subsection{Hamiltonian-like vector fields on almost-symplectic
manifolds} 
An almost-symplectic manifold $(M,\s)$ is a generalization of a
symplectic manifold, in which a smooth manifold $M$ of dimension
$2n$, with $n\ge2$, is equipped with a non-degenerate but not closed 2-form $\s$.
Hamiltonian-like systems on almost-symplectic manifolds arise, for
instance, in nonholonomic mechanics \cite{bates-sn,koiller};
moreover, their study might be an intermediate step to the study of
the more general case of (generalized) Hamiltonian-like systems on
twisted almost-Dirac manifolds considered in \cite{BlvdS,SW,zung2}. 
There is however a crucial difference from the symplectic context. 

A Hamiltonian vector field $X$ on a symplectic manifold $(M,\s)$ is
commonly defined as a vector field such that the contraction $i_X\s$
is exact. By the closedness of $\s$, this definition is equivalent to
the vanishing of the Lie derivative $L_X\s=i_Xd\s+d(i_X\s)$, that is,
to the fact that $X$ is a symmetry of the symplectic form. This fact
has far reaching consequences, being the reason of the Lie algebra
structure of the (locally) Hamiltonian vector fields, which is
typical of the symplectic world and dominates, e.g., the field of
complete integrability. 

In the almost-symplectic setting this link is broken. One may still
retain as definition of ``being Hamiltonian'' the exactness of the
contraction $i_X\s$, as e.g. in \cite{cannas}, but then, since
$d\s\not=0$, these vector fields need not be symmetries of $\s$. This
destroys the Lie algebra structure of the (locally) Hamiltonian vector
fields and is related in an obvious way to the failure of the Jacobi identity
of the almost-Poisson bracket induced by an almost-symplectic
structure. 

Nevertheless there may exist, in an almost-symplectic manifold,
special vector fields which have both properties $d(i_X\s)=0$ and
$L_X\s=0$. These vector fields, that closely resemble the Hamiltonian
vector fields of the standard symplectic case, were called
``strongly-Hamiltonian'' in~\cite{FS1} and ``Hamiltonian'' in
\cite{vaisman}. In view of the fact that they have the two
qualifying properties of the Hamiltonian vector fields of the
symplectic case, here we prefer to call them {\it
fully-Hamiltonian}. Correspondingly, in order to avoid ambiguities, we
will call {\it partially-Hamiltonian} those vector fields $X$
for which $i_X\s$ is exact. We will occasionally
call ``standard Hamiltonian'' those vector fields which are
Hamiltonian with respect to a symplectic structure.

Fully-Hamiltonian vector fields were identified in \cite{FS1} (even
though they had already appeared in \cite{BlvdS} as infinitesimal
generators of group actions in the framework of almost-Dirac
structures). Subsequently, they have
been considered in \cite{vaisman,SS1,fok}. 
In particular, Vaisman \cite{vaisman} points out that 
fully-Hamiltonian vector fields are rather special and that, since on
a given almost-symplectic manifold there might exist very few or
even none of them which are non-zero, it is of interest to provide
examples of the existence and non-existence of these vector fields and
investigate their structure, and does it. The present paper provides
contributions in this direction.

\subsection{Aim of the paper}
Specifically, we study the structure of fully-Hamiltonian vector
fields on a particular class of almost-symplectic manifolds, which
was considered in \cite{FS1} in connection with a generalization of
the known patterns of integrability from the symplectic to the
almost-symplectic context. The almost-symplectic manifolds we consider are 
equipped with a fibration whose fibers are Lagrangian tori
(namely, $\s$ vanishes on them) and has the
additional property that the partially-Hamiltonian vector field of any
basic (namely, constant on the fibers) function is fully-Hamiltonian.
We will call {\it fully-Lagrangian toric fibration} a fibration of
this type. The {\it semi-global} structure of these fibrations (namely,
their structure in a saturated neighbourhood of each
fiber) has been studied in \cite{FS1}. In particular, it is there
proven that only cases with $n:=\frac12\dim M\ge3$ are interesting,
because if $n=2$ then $\s$ is necessarily symplectic. 

On any such fibered manifold there are two distinct classes of 
fully-Hamiltonian vector fields: those which are vertical (namely, 
tangent to the fibers of the fibration) and those which are not
vertical, and we are mostly interested to their
semi-global structure that we will investigate through the introduction of
action-angle coordinates.

{\it Vertical} fully-Hamiltonian vector fields have already been considered
in \cite{FS1}. They are exactly the partially-Hamiltonian vector
fields whose Hamiltonians are basic functions and are integrable, in
the sense that their dynamics is conjugate to a linear flow on
$\tenne$, with $n=\frac12\dim M$. Moreover, semi-globally,
they are Hamiltonian (and completely integrable, in the sense
that they possess $n$ almost everywhere independent
Poisson commuting first integrals) with respect to a
semi-global {\it symplectic} structure defined by the
almost-symplectic structure. Thus, the semi-global structure of these
fully-Hamiltonian vector field is not particularly interesting.
Globally, things could be different: if the fibration is
topologically non trivial, then there might not exist a global
symplectic structure in $M$ with respect to which all vertical
fully-Hamiltonian vector fields are Hamiltonian. However, we leave
this global question for future investigations.

The semi-global structure of {\it non-vertical} fully-Hamiltonian vector
fields is much more interesting, and is the main topic of this work.
We will make certain hypotheses on the fully-Lagrangian toric fibration
(it must be produced by a free ``fully-Hamiltonian'' action of $\toro^n$)
and on the fully-Hamiltonian vector fields we look for (which involve
certain genericity properties of the Fourier spectrum of their
Hamiltonians). In short, we will show that, under the considered
hypotheses:
\begin{itemize}
\item[i.] If $n=3$, then there are no non-vertical fully-Hamiltonian
vector fields (Theorem \ref{thm1}).
\item[ii.] If $n=4$, then, semi-globally, every non-vertical fully-Hamiltonian 
vector field is Hamiltonian (and completely
integrable) with respect to some semi-globally defined 
symplectic structure (Theorem \ref{thm3}).
\item[iii.] For any $n\ge4$, every non-vertical fully-Hamiltonian 
vector field is equivariant under a $\toro^k$-sub-action of the
$\tenne$-action that produces the fibration, for some $3\le k\le n$,
and reduces to a family of standard Hamiltonian vector fields on
reduced {\it symplectic} manifolds of dimension $2n-2k$ (Theorem
\ref{thm2}). The reduction process uses an almost-symplectic version
of the standard Marsden-Meyer-Weinstein symplectic reduction
procedure due to \cite{BlvdS,vaisman}. 
\item[iv.] If $n\ge5$, non-vertical fully-Hamiltonian vector fields
need not be integrable.
\end{itemize}
We will provide examples of all these cases.

In the current terminology, particularly in the field
of nonholonomic mechanics (see e.g. \cite{ehlers, fedorov, luis,
bolsinov, balseiro, borisov}), a result like iii. would be considered
a {\it Hamiltonization} result. In the terminology that we adopt
in this work, which focuses on the difference between symplectic and
almost-symplectic, we rather regard it as a {\it symplectization}
result.

\begin{remarks} 
\begin{itemize} 

\item[(i.)] Motivations for this work come also from dynamics, specifically,
from the
study of almost-symplectic perturbations of standard completely integrable
Hamiltonian systems, a topic of which very little is known. For certain
types of perturbations, such systems can be seen as defined in a
fully-Lagrangian toric fibration (created by the standard integrable
Hamiltonian system). Our results indicate that the case of
fully-Hamiltonian and partially-Hamiltonian perturbations could be
very different, since the former could be reduced to that of
a family of standard Hamiltonian perturbations of standard integrable
Hamiltonian systems. 

\item[(ii.)] A class of almost-symplectic systems which has received a
certain attention is that of conformally symplectic systems (see e.g.
\cite{vaisman-CS,liverani,banyaga,maro}). 
However almost-symplectic manifolds that carry fully-Lagrangian toric fibrations 
are not conformally symplectic symplectic.
\end{itemize}
\end{remarks}
 
\subsection{Organization of the paper} 

Sections 2 and 3 introduce a number of topics which are necessary for
the present study. In particular, Section \ref{s:2} has two purposes.
The first is to recall---and here and there to develop---some basic
facts about almost-symplectic systems and almost-symplectic Lie group
actions and reduction. The second is to prove a general constraint on
fully-Hamiltonian vector fields on any almost-symplectic manifold:
roughly speaking, their Hamiltonian functions are locally independent
of at least three coordinates (Propositions \ref{p:upperbound} and
\ref{p:upperbound-2}). To a certain extent, this fact is the ultimate
reason for our symplectization result.

In Section \ref{s:3}, recalling and extending \cite{FS1}, we
introduce the ``fully-Lagrangian toric fibrations'' and we add a
characterization of these structures in terms of invariance under a
fully-Hamiltonian torus action, which is needed in view of the
symplectization procedure (Proposition \ref{p:converse}).

In Section \ref{s:4} we first introduce the genericity assumptions that
enter our symplectization results and then state them
(Theorems \ref{thm1}, \ref{thm3} and \ref{thm2}). The genericity
assumptions require a Fourier series on torus bundles, that we
introduce following previous ideas from \cite{FF-hpt,FF-cg,FFS,FF-Enc}.
Theorems \ref{thm1} and \ref{thm2} are then proved in Section
\ref{s:proofs}.

Section \ref{s:esempi} describes and comments the ``reconstruction''
procedure needed to deduce the dynamical properties of non-vertical
fully-Hamiltonian vector fields from those of the reduced vector
fields.

Unless said differently, all objects are assumed to be smooth
and all vector fields are assumed to be complete. ``Fibration''
means ``locally trivial fibration''. The flow of a vector field $X$
will be denoted $\P^X$, with time-$t$ maps $\P^X_t$. Of course, all
angles and their sums have to be understood $\mathrm{mod}\, 2\pi$.

\section{Hamiltonian-like systems in almost-symplectic manifolds }
\label{s:2}

\subsection{Partially- and fully-Hamiltonian vector fields}
An {\it almost-symplectic} manifold $(M,\s)$ is an even-dimensional
manifold $M$ equipped with a non-degenerate 2-form~$\s$. Note that we
do not require $\s$ to be closed, and thus regard symplectic manifolds as 
particular almost-symplectic manifolds. 

If $\dim M=2n$, then we say that $n$ is the {\it number of degrees of
freedom} of $(M,\s)$. The case of interest is $n\ge 2$, because when
$n=1$ every 2-form is closed and $(M,\s)$ is symplectic. Nevertheless,
since as a result of reductions under group actions we will reach
cases with $n=1$ or even $0$, we include the case $n=1$ in all our
definitions relative to almost-symplectic manifolds, and make the
convention that if $M$ is zero-dimensional then $(M,\s)$, with
$\s=0$, is a symplectic manifold.

\begin{definition}
\label{d:Def1}
A vector field $X$ on an almost-symplectic manifold $(M,\s)$ is:
\begin{itemize}
\item[i.] {\it Partially-Hamiltonian} if $i_X\sigma$ is exact. Thus
$i_X\sigma= -df$ for some function $f\in C^{\infty}(M)$ that we call a
{\it Hamiltonian} of $X$, and we write $X^\s_f$ for $X$. 
\item[ii.] {\it Fully-Hamiltonian} if, besides being
partially-Hamiltonian, it is a symmetry of $\sigma$, namely
$L_X\sigma = 0$.
\end{itemize}
If $X^\s_f$ is fully-Hamiltonian, then we say that 
$f$ is a {\it fully-Hamiltonian function on $(M,\s)$}. The
set of all fully-Hamiltonian functions on $(M,\s)$ is denoted
$\cH^\s_F(M)$.
\end{definition}

As already noticed, $L_X\s=i_X(d\s)+d(i_X\s)$ and the
full-Hamiltonianity of a partially-Hamiltonian vector field is
equivalent to
\begin{equation}
\label{ixdsigma}
   i_Xd\sigma= 0 \,.
\end{equation}
Clearly, $L_{X^\s_f}f=0$ for any function $f$, so every
partially-Hamiltonian vector field has its own Hamiltonian as a first
integral.

We will occasionally use also local versions of these objects:

\begin{definition}
\label{d:Def2}
A vector field $X$ on an almost-symplectic manifold $(M,\s)$ is {\it
locally partially-Hamiltonian} if $i_X\s$ is closed, and {\it locally
fully-Hamiltonian} if, moreover, \hbox{$i_Xd\s=0$}. 
\end{definition}

At an algebraic level, the reason for considering fully-Hamiltonian
vector fields is the following. The almost-symplectic form $\sigma$
induces an almost-Poisson bracket on smooth functions of $M$, which
is defined by
\begin{equation}
\label{eq:aP}
   \{f,g\}^\s := - \sigma(X^\s_f,X^\s_g) \qquad \forall\,f,g\in C^\infty(M) \,.
\end{equation}
As in the symplectic case, $\{f,g\}^\s=-L_{X_f^\s}g$. However,
because of the non-closedness of $\sigma$, this bracket does not
satisfy the Jacobi identity \cite{cannas}. Therefore, it does not make
$C^\infty(M)$ a Lie algebra and does not induce an
\hbox{(anti-)}homomorphism between functions and
partially-Hamiltonian vector fields on $M$. However, there is an
analogue of this homomorphism for the restriction of the bracket to
$\cH^\s_F(M)$. This is a consequence of the following fact (from
\cite{FS1}, where its statement contains an obvious flaw that we
fix here):

\begin{proposition}
\label{l:homomorphism} 
Let $Y$ and $Z$ be two vector fields on an almost-symplectic manifold
$(M,\sigma)$. If $Y$ is locally partially-Hamiltonian and $Z$ is a
symmetry of $\s$ then $[Y,Z]$ is partially-Hamiltonian with
Hamiltonian $-\sigma(Y,Z)$:
$$
    [Y,Z] = - X^\s_{\sigma(Y,Z)} \,.
$$
\end{proposition}

\begin{proof} Since $d(i_Y\s)=0$ and $L_Z\s=0$,
$ d(\sigma(Y,Z)) = d(i_Zi_Y\s) =
L_Z(i_Y\sigma) - i_Z\, d(i_Y\sigma) = i_Y(L_Z\sigma)+i_{[Z,Y]}\sigma
= i_{[Z,Y]}\sigma = - i_{[Y,Z]}\sigma$.
\end{proof}

\begin{corollary}
\label{c:anti-hom}
If two vector fields $Y$ and $Z$ are locally fully-Hamiltonian, then
their Lie bracket $[Y,Z]$ is fully-Hamiltonian.
\end{corollary}

\begin{proof} If $Y$ and $Z$ are locally fully-Hamiltonian,
then $[Y,Z]=-X_{\s(Y,Z)}$ is partially-Hamiltonian and
$L_{[Y,Z]}\s= (L_Y\circ L_Z -L_Z\circ L_Y)\s=0$.
\end{proof}

Therefore, the set of all locally fully-Hamiltonian vector fields
is a Lie subalgebra of the Lie algebra of all vector fields on $M$.
In particular, for all fully-Hamiltonian functions $f,g\in\cC^\infty(M)$,
$$
   [X^\s_f,X^\s_g] = - X^\s_{\{f,g\}^\s} \,,
$$
so that, equipped with the bracket \for{eq:aP}, $\cH^\s_F(M)$ is a Lie algebra
and the map $f\mapsto X^\s_f$ is an anti-homomorphism between it and
the Lie algebra of fully-Hamiltonian vector fields.

In the class of partially-Hamiltonian vector fields, those vector
fields which are fully-Hamiltonian have also special dynamical
properties. For instance, partially-Hamiltonian vector fields need not
preserve the volume form $\s^n$: an example is
$X=\der{}{x_3}$ on $M=\reali^4\setminus\{0\}$ with
$\s=x_3dx_1\wedge dx_2+dx_3\wedge dx_4$. Instead, since
$L_X\s^n=n\s^{n-1}\wedge L_X\s$:

\begin{proposition}
\label{p:ConsVol} Every locally fully-Hamiltonian vector field on an
almost-symplectic manifold $(M,\s)$ preserves the volume $\s^n$.
\end{proposition}

\begin{remark} 
There are many studies on the Hamilton-Jacobi theory outside the
standard symplectic framework, see e.g. \cite{CGMMLM,BMMdDP,EDLLSZ}.
However, as pointed out in \cite{BFS} in connection with the almost-symplectic formulation of
non-holonomic systems \cite{bates-sn}, the absence of the Lie algebra
(anti-)homomorphism between partially Hamiltonian
functions and partially Hamiltonian vector implies that, for an
almost-symplectic system, the existence of a
complete integral of the Hamilton-Jacobi equation does not imply,
without further hypotheses related to full-Hamitonianity, 
the integrability of the system.
\end{remark}

\subsection{An estimate on the number of (independent) locally
fully-Hamiltonian vector fields}
\label{ss:Estimate} 
As remarked by \cite{vaisman}, locally fully-Hamiltonian vector fields
on $(M,\s)$ define a generalized distribution $\cD^\s_F$ on $M$,
whose fiber at each point $m\in M$ is the set of the values in $m$ of
all the locally fully-Hamiltonian vector fields. This set is clearly
a subspace of $T_mM$. As proven in \cite{vaisman}, this generalized
distribution is integrable and endows $M$ with a Dirac structure.
Here, we add the following simple quantitative estimate, which will
play a key role in the subsequent analysis. 

\begin{proposition}
\label{p:upperbound} $\cD_F^\s$ has rank $\le 2n-3$ at each point
$m\in M$ at which $d\s(m)\not=0$.
\end{proposition}

\begin{proof} 
By \for{ixdsigma}, all vectors $u\in T_mM$ which belong to the
fiber of $\cD_F^\s$ over $m$ annihilate the 3-form $d\sigma$ at
$m$, that is, $i_ud\s(m)=0$. The statement now follows from the fact
that the kernel of a $p$-form $\eta$ has codimension $\ge p$ at each
point at which it is nonzero. [If $\eta(m)\not=0$, then there
exist $p$ linearly independent vectors $u_1,\ldots,u_p \in T_mM$ such
that $\eta(m)(u_1,\ldots,u_p)\not=0$. Therefore
$i_{u_j}\eta(m)\not=0$ for all $j=1,\ldots,p$ and every set of
linearly independent vectors that annihilate $\eta(m)$ has
cardinality $\le 2n-p$.] 
\end{proof}

This implies that, locally, in a neighbourhood of each point at which
$d\s$ is non-zero, there are no more than $2n-3$ linearly independent
(germs of) locally fully-Hamiltonian vector fields. 

Under suitable stronger conditions it is possible to
say more about the structure of the (local) Hamiltonian of the locally
fully-Hamiltonian vector fields. For instance:

\begin{proposition}
\label{p:upperbound-2}
Assume that $d\s(m)\not=0$ at a point $m\in M$ and that in a
neighbourhood of $m$ the generalized distribution $\cD_F^\s$ is
generated by a finite set of $2n-s$ everywhere linearly independent
locally fully-Hamiltonian vector fields. 

Then, $s\ge3$ and there exists a coordinate system in a neighbourhood
$V$ of $m$ such that the restriction to $V\cap U$ of any local
Hamiltonian $h\in C^\infty(U)$, with $U$ a neighbourhood of $m$, of a
locally fully-Hamiltonian vector field is independent of the last $s$
coordinates.

\end{proposition}

\begin{proof}
That $s\ge3$ follows from Proposition \ref{p:upperbound} because, under
the stated hypotheses, $\cD_F^\s$ has rank $2n-s$ at $m$.

Choose a set of locally fully-Hamiltonian vector fields
$X_1,\ldots,X_r$ which generate $\cD_F^\s$ in a neighbourhood $V_0$ of
$m$. In a neighbourhood $V_1\subseteq V_0$ of $m$ these vector fields
have local Hamiltonians $f_1,\ldots,f_r\in C^\infty(V_1)$,
$X_j=X_{f_j}^\s$. By the non-degeneracy of $\s$, $f_1,\ldots,f_r$
are everywhere functionally independent in $V_1$ and can be completed
to a coordinate system $(f_1,\ldots,f_r,g_{1},\ldots,g_s)$ in a
neighbourhood $V\subseteq V_1$ of $m$. Since
$X^\s_{f_1},\ldots,\ldots,X^\s_{g_s}$ are everywhere linearly
independent in $V$ and the $X^\s_{f_j}$'s generate $\cD_F^\s$ at each
point of $V_2$, the $X^\s_{g_j}$ do not belong to $\cD^\s_F$ at any
point of $V_2$. 

Consider now a locally partially-Hamiltonian vector field $Y$, with 
a local Hamiltonian $h$ in a neighbourhood $U$ of $m$. In $V\cap U$,
$$
  Y = \sum_{j=1}^r\der h {f_j}X^\s_{f_j} + 
        \sum_{j=1}^s\der h {g_j} X^\s_{g_j} \,.
$$
Thus, $Y|_{V\cap U}$ is a section of $\cD_F^\s$, namely locally
fully-Hamiltonian, if and only if all $\der h{g_j}=0$.
\end{proof}

Even though we will not actually use this result in the subsequent
analysis, it should be kept in mind as an underlying motivation
for it.

We note that the upper bound $2n-3$ on the maximum number of independent
fully-Hamiltonian functions given in Proposition
\ref{p:upperbound} is {\it de facto} met in all examples in
\cite{vaisman}, and so is the existence of a coordinate system as in
Proposition \ref{p:upperbound-2}. 

\subsection{Almost-symplectic diffeomorphisms, group actions, and
reduction} 
\label{ss:ASdiffeos}

We conclude this section with some (rather obvious) remarks on the
almost-symplectic analogues of symplectomorphisms, of symplectic and
Hamiltonian group actions, and of the reduction procedure which will
be needed in the following.

\begin{definition}
An {\it almost-symplectic diffeomorphism} between two
almost-symplectic manifolds $(M,\s)$ and $(M',\s')$ is a
diffeomorphism $\rho:M\to M'$ such that $\rho_*\s=\s'$.
\end{definition}

\begin{proposition} 
\label{p:ASD-1}
If $\rho$ is an almost-symplectic diffeomorphism of $(M,\s)$ onto
$(M',\s')$, then:
\begin{itemize}
\item[i.] $\rho_*X^\s_f=X^{\s'}_{\rho_*f}$ for all
$f\in\cC^\infty(M)$.
\item[ii.] If $f\in\cH_F^\s(M)$, then $\rho_*f\in\cH_F^{\s'}(M')$.
\end{itemize}
\end{proposition}

\begin{proof} Write $X$ for $X_f^\s$. Then $i_{\rho_*X}(\rho_*\s)=
\rho_*(i_X\s) = -\rho_*df = - d(\rho_*f)$ and 
$L_{\rho_*X}(\rho_*\s)=\rho_*(L_X\s)=0$. 
\end{proof}

\begin{proposition} 
\label{p:ASD-2}
If $X$ is a locally fully-Hamiltonian vector field on an
almost-symplectic manifold $(M,\s)$, then, for each $t\in\reali$,
$\P^{X}_t$ is an almost-symplectic diffeomorphism of $(M,\s)$ onto
itself.
\end{proposition}

\begin{proof} From $\frac d{dt}\big(\P^X_t)_*\s = (\P^X_t)_*L_X\s=0$
it follows $(\P^X_t)_*\s=(\P^X_0)_*\s= \s$. \end{proof}

We consider now an action $\Psi$ of a Lie group $G$ on an
almost-symplectic manifold $(M,\s)$. As usual, we denote $\mathfrak
g$ the Lie algebra of $G$ and $\exp$ its exponential map. For any
$\xi\in\mathfrak g$, we denote $Y_\xi$ the infinitesimal generator
$\frac d{dt}\Psi_{\exp(t\xi)}\big|_{t=0}$. 

\begin{definition}
\label{d:Def4}
The action $\Psi$ of $G$ on $M$ is:

\begin{itemize}

\item[i.] {\it Almost-symplectic} if $(\Psi_g)_*\s=\s$ for all
$g\in G$.

\item[ii.] {\it (Locally) partially-Hamiltonian} if, for all
$\xi\in\mathfrak g$, $Y_\xi$ is (locally) partially-Hamiltonian.

\item[iii.] {\it (Locally) fully-Hamiltonian} if, for all
$\xi\in\mathfrak g$, $Y_\xi$ is (locally) 
fully-Hamiltonian.
\end{itemize}
\end{definition}

Fully-Hamiltonian actions are called {\it
strongly-Hamiltonian} in \cite{BlvdS} and {\it Hamiltonian} in
\cite{vaisman}.

\begin{proposition}
\label{p:GA-1}
\begin{itemize}
\item[i.] Any almost-symplectic and (locally) partially-Hamiltonian
action is (locally) fully-Hamiltonian.
\item[ii.] Assume that the exponential map of $G$ is surjective. Then
any locally fully-Hamiltonian action is almost-symplectic.
\end{itemize}
\end{proposition}

\begin{proof} (i.) Assume that $\Psi$ is almost-symplectic. Then, 
$\forall \xi\in\mathfrak g$, $(\Psi_{\exp(t\xi)})_*\s=\s$ for all $t$
and hence, since $\Psi_{\exp(t\xi)} = \P^{Y_\xi}_t$, $L_{Y_\xi}\s=0$.
It follows that, if $\Psi$ is (locally) partially-Hamiltonian, then it is
(locally) fully-Hamiltonian. 

(ii.) If the action is locally fully-Hamiltonian, then all
infinitesimal generators $Y_\xi$, $\xi\in\mathfrak g$, are locally
fully-Hamiltonian. Therefore, by Proposition \ref{p:ASD-2},
$\Psi_{\exp(\xi)}$ is an almost-symplectic diffeomorphism
$\forall\xi\in\mathfrak g$ and, if $\exp$ is surjective, so is
$\Psi_g$ for all $g\in G$. \end{proof}

A partially-Hamiltonian action has a momentum map
$$
   J:M\to\mathfrak g^*
$$
which is defined by the fact that, for each $\xi\in\mathfrak g$, the
function $J_\xi:=\langle J,\xi\rangle$ (with
$\langle \ ,\ \rangle$ the duality $\mathfrak g$-$\mathfrak g^*$) satisfies
$
   dJ_\xi=-i_{Y_\xi}\s \,,
$
and hence $Y_\xi=X_{J_\xi}^\s$. If the action is fully-Hamiltonian,
then so is its momentum map, in the sense that all functions $J_\xi$
are fully-Hamiltonian. Moreover (see \cite{BlvdS,FS1}):

\begin{proposition} 
\label{p:mm1} 
Consider a partially-Hamiltonian $G$-action $\Psi$ on an
almost-symplectic manifold $(M,\s)$ with momentum map~$J$. Let
$f\in\cC^\infty(M)$ be a $G$-invariant function. Then:
\begin{itemize}
\item[i.] $J$ is constant along the flow of~$X_f^\s$.
\item[ii.] If either the action $\Psi$ or the function $f$ is
fully-Hamiltonian, then $X_f^\s$ is $G$-invariant.
\end{itemize}
\end{proposition}

\begin{proof} (i.) If $f$ is $G$-invariant then
$L_{X_f^\s}J_\xi=-i_{X_f^\s}i_{Y_\xi}\s=-L_{Y_\xi}f$ vanishes for all
$\xi\in\mathfrak g$.

(ii.) Fix $\xi\in\mathfrak g$. Since at least one among the two
partially-Hamiltonian vector fields $X_f^\s$ and $Y_\xi=X^\s_{J_\xi}$
is fully-Hamiltonian, Proposition \ref{l:homomorphism} gives
$[X_f^\s,Y_\xi]=-X_{\{f,J_\xi\}^\s}^\s$. But $\{f,J_\xi\}^\s =
L_{X_f^\s}J_\xi$ which, by item i., vanishes.
\end{proof}

As usual, we say that a momentum map $J$ of $\Psi$ is {\it
equivariant} if it is equivariant with respect to the action $\Psi$
on $M$ and to the coadjoint action of $G$ on $\mathfrak g^*$. The
equivariance of $J$ implies that, for each $c\in\mathfrak g^*$, $G_c$
(the isotropy subgroup of $c$ relative to the coadjoint action of
$G$) acts on $J^{-1}(c)$. 

If $J$ is equivariant, $c\in\mathfrak g^*$ is a regular value of
$J$ and the action of $G_c$ on $J^{-1}(c)$ is free and proper, then
$$
  M_c:=J^{-1}(c)/G_c
$$
is a smooth manifold and the canonical projection $p_c:J^{-1}(c) \to
M_c$ is a smooth submersion. Furthermore, given a $G$-invariant
function  $f\in\cC^\infty(M)$, there exists a unique function
$f_c\in\cC^\infty(M_c)$ such that
$$
  p_c^*f_c=\imm_c^*f 
$$
where $\imm_c:J^{-1}(c) \hookrightarrow M$ is the canonical inclusion.
Slightly extending previous results by \cite{BlvdS,vaisman} (which
state only item i.) we have
the following almost-symplectic version of the well known
Marsden-Meyer-Weinstein symplectic reduction theorem:

\begin{proposition} 
\label{p:vaisman}
Assume that $\Psi$ is almost-symplectic and
fully-Hamiltonian, that $J$ is equivariant,
that $c$ is a regular value of $J$ and that  $G_c$ acts freely and
properly on $J^{-1}(c)$. Then:
\begin{itemize}
\item[i.] There exists an almost-symplectic form $\s_c$ on $M_c$ such
that $p_c^*\s_c=\imm_c^*\s$. 
\item[ii.] If $f\in\cC^\infty(M)$ is $G$-invariant, then
$X_{f_c}^{\s_c}$ is $p_c$-related to the restriction of $X_{f}^{\s}$
to $J^{-1}(c)$  (namely, $T_mp_c \cdot X_f^\s(m) =
X_{f_c}^{\s_c}(p_c(m))$ for all $m\in \imm(J^{-1}(c))$).
\item[iii.] If $f\in\cH_F^\s(M)$, then $f_c\in\cH_F^{\s_c}(M_c)$.
\end{itemize}
\end{proposition}

\begin{proof} (i.) See \cite{vaisman}. (ii.) This does not need the
closedness of $\s$ and its proof goes as in the symplectic case, see
e.g. Theorem 14.6 in \cite{libmarle}. (iii.) Using i., 
$p_c^*(L_{X^{\s_c}_{f_c}}\,\s_c) =  L_{X^\s_f|_{\imm(J^{-1}(c))}}\,
\imm_c^*\s = \imm_c^*(L_{X^\s_f}\s) =0$.~\end{proof}

We will say that the $(M_c,\s_c)$'s, $f_c$'s and $(M_c,\s_c,f_c)$'s
are, respectively, the {\it almost-symplectic reduced spaces}, {\it
reduced Hamiltonians} and {\it almost-symplectic reduced systems}
induced by $(M,\s,f)$.

\begin{remarks} 
\begin{itemize}
\item[(i.)] It is possible that, for some or all $c$,
$d\s_c=0$ and the reduced almost-symplectic space $(M_c,\s_c)$ is in
fact symplectic. 

\item[(ii.)] The reduced manifolds $M_c$ might not be connected even if $M$
is connected. This is a delicate point, which has been studied in the
symplectic case (see \cite{Atiyah,GuilleminSternberg,Kirwan}); for
instance, in that case, if $G=\toro^k$ then the connectedness is
guaranteed if $M$ is compact, but not in general).
\end{itemize}
\end{remarks}

\section{Fully-Lagrangian toric fibrations on almost-symplectic
manifolds}
\label{s:3}

In this section we introduce the class of
fully-Hamiltonian vector fields that we study in the rest of the
paper; essentially, we follow, and globalize, \cite{FS1}.

\subsection{Fully-Lagrangian toric fibrations}
\label{ss:SHCIS}
A submanifold $N$ of an almost-symplectic manifold $(M,\s)$ is said
to be {\it isotropic} if the restriction of $\s$ to it vanishes and
{\it Lagrangian} if it is isotropic and
$\mathrm{dim}(N)=\frac12\mathrm{dim}(M)$.  A submanifold $N$ of $M$
is Lagrangian if and only if, for each $x\in N$, $T_xN$ equals its
own $\s$-orthogonal $(T_xN)^\s:=\{u\in T_xM: \s(u,v)=0\ \forall \,
v\in T_xN\}$.  Clearly, if a fibration $\pi:M\to B$ has Lagrangian
fibers then $\dim M=2\dim B$. 

Given a fibration $\pi:M\to B$, a $p$-form ($p\ge0$) on the total
space $M$ is said to be {\it basic} if it is the pull-back under
$\pi$ of a $p$-form on the base $B$. Basic functions ($p=0$) are the
functions on $M$ which are constant on the components of the fibers
of $\pi$. A vector field on $M$ is said to be {\it vertical} if it is a section
of $\ker(T\pi)$, namely, it is tangent to the fibers of $\pi$.

\begin{definition}
\label{d:FLF} 
Let $(M,\s)$ be an almost-symplectic manifold. A {\it fully-Lagrangian
toric fibration} of $(M,\s)$ is a fibration $\pi : M \to B$ such that
\begin{itemize}
\item[i.] The fibers of $\pi$ are compact and connected Lagrangian
submanifolds of $M$.
\item[ii.] All basic functions are fully-Hamiltonian.
\end{itemize}
\end{definition}

An example of a fully-Lagrangian fibration of an almost-symplectic
manifold $(M,\s)$ is provided by a submersion
$$
    \Pi=(\Pi_1,\ldots,\Pi_n):M\to\Pi(M)\subseteq\reali^n
$$
with compact and connected fibers whose components $\Pi_1,\ldots,\Pi_n$ are
fully-Hamiltonian functions and are pairwise in involution, namely
$$
   \{\Pi_i,\Pi_j\}^\s=0 \qquad \forall\,i,j \,.
$$
Indeed, the fibers of any such $\Pi$ are obviously Lagrangian (their
tangent spaces are spanned by the $X^\s_{\Pi_i}$'s because
$L_{\Pi_i}\Pi_j=-\{\Pi_i,\Pi_j\}^\s=0$, and 
$\s(X^\s_{\Pi_i},X^\s_{\Pi_j})=-{\{\Pi_i,\Pi_j\}^\s}=0$) and all basic
functions are fully-Hamiltonian (if $F:\Pi(M)\to\reali$, then
$X^\s_{F\circ\Pi} = \sum_i(\der F{\Pi_i}\circ\Pi) X^\s_{\Pi_i}$ and so
$i_{X^\s_{F\circ\Pi}}d\s=0$). This is the case considered in
\cite{FS1}.

As it turns out, a submersion $\Pi$ with this properties is a
``semi-global'' (that is, in a saturated neighbourhood of any fiber)
model for any fully-Lagrangian toric fibration $\pi:M\to B$. In fact,
if $\cB=(\cB_1,\ldots,\cB_n)$ is any set of local coordinates on the
base $B$, then the components of the submersion $\Pi:=\cB\circ\pi$ are
basic and hence fully-Hamiltonian functions, while the fact that they
are pairwise in involution follows from the following fact (which is
a known property of Lagrangian fibrations of symplectic manifolds and
extends to the almost-symplectic context because it does not rely on
the closedness of $\s$; see also Proposition 7 of \cite{FS1}):

\begin{lemma}\label{p:basicfunctions-1}
Let $\pi:M\to B$ be a fully-Lagrangian toric fibration on an
almost-symplectic manifold $(M,\s)$. A function $f$ on $M$ is basic
if and only if any of the following two equivalent conditions is
satisfied:
\begin{itemize}
\item[i.] $X_f^\s$ is vertical.
\item[ii.] $f$ is in involution with all basic functions.
\end{itemize}
\end{lemma}

\begin{proof} 
(i.) $f$ is basic $\Leftrightarrow$ $f$ is a first integral of $\pi$ 
$\Leftrightarrow$ $\s(Y,X^\s_f) = L_Yf= 0$ for all vertical vector
fields $Y$ $\Leftrightarrow$ $X^\s_f$ is
$\s$-orthogonal to the fibers of $\pi$ $\Leftrightarrow$ $X^\s_f$ is
tangent to the fibers of $\pi$, given that they are Lagrangian
submanifolds.

(ii.) If $f$ and $g$ are basic then, by i., $X_f^\s$ and $X^\s_g$ are
vertical and, since the fibers are Lagrangian, $\{f,g\}^\s =
-\s(X_f^\s,X_g^\s) = 0$. Conversely, if $\{f,g\}^\s=0$ with a basic
function $g$, then $L_{X_f^\s}g\not=0$ and $X_f$ is not tangent to the
fibers of $g$, and thus neither to those of $\pi$, which are contained in them.
\end{proof}

The requirement of the full-Hamiltonianity of the basic functions
of a fully-Lagrangian toric fibration $\pi:M\to B$
implies that its fibers have a set of everywhere linearly
independent tangent vector fields (e.g., for any set
of local coordinates $\cB_1,\ldots,\cB_n$ on $B$,
the vector fields $X_{\cB_i\circ\pi}^\s$'s) which are pairwise in
involution ($[X_{\cB_i\circ\pi}^\s, X_{\cB_j\circ\pi}^\s] = 
X_{ \{\cB_i\circ\pi,\cB_j\circ\pi\}^\s} =0$) and thus, by a well
known argument, are diffeomorphic to $\tenne$. Even more,
an almost-symplectic version of the Liouville-Arnold theorem is
valid, which describes the semi-global almost-symplectic structure of
such fibrations:

\begin{proposition}
\label{p:CI} {\rm \cite{FS1} } Let $\pi:M\to B$ be a
fully-Lagrangian toric fibration of an almost-symplectic manifold
$(M,\s)$ with $\dim(M)=2n$ and $n\ge1$. 
\begin{itemize}
\item[i.] The fibers of $\pi$ are diffeomorphic to $\tenne$.
\item[ii.] Each fiber of $\pi$ has a neighbourhood $V$ equipped with
a diffeomorphism $(a,\alpha): V\to\cA\times\tenne \subseteq 
\renne\times\tenne$, such that the fibers of $\pi$ coincide with the
sets $a=\const$ and the local representative $\saa:=(a,\a)_*\s$ of
$\s$ has the form
\begin{equation}
\label{eq:sigma}
  \saa(a,\a) = \sum_i da_i\wedge d\alpha_i
  + \frac12 \sum_{ij} A_{ij}(a) da_i\wedge da_j
\end{equation}
where $A(a)$ is an $n\times n$ skew-symmetric matrix that depends
smoothly on $a$.
\end{itemize}
\end{proposition}

The trivializations $(a,\alpha)$ provided by Proposition \ref{p:CI}
will be called (local) {\it action-angle coordinates} of the
fibration $\pi$.  Written in these coordinates, the
partially-Hamiltonian vector field $X_f^\s$ of a function $f$ has
components
\begin{equation}\label{eq:HVF}
  \big(X_f^\s\big)_{a_i} = - \der \cF {\alpha_i} \,,\qquad
  \big(X_f^\s\big)_{\alpha_i}  = \der \cF{a_i} + \sum_jA_{ij}\der \cF{\alpha_j} 
  \,, \qquad i=1,\ldots,n \,,
\end{equation}
where $\cF$ is the local representative of $f$. Moreover, if $\cF_1$ and
$\cF_2$ are the local representatives of two functions $f_1$ and $f_2$ in an
action-angle chart, then the local representative of the
almost-Poisson bracket $\{f_1,f_2\}^\s$ in that chart turns out to be
\begin{equation}\label{eq:APB}
  \{\cF_1,\cF_2\}^{\saa} :=
  \sum_i\Big( \der {\cF_1}{a_i} \der{\cF_2}{\alpha_i} - 
                    \der{\cF_1}{\alpha_i} \der{\cF_2}{a_i} \Big)
  + 
  \sum_{ij} A_{ij}\der{\cF_1}{\alpha_i} \der{\cF_2}{\alpha_j} \,.
\end{equation}

As in the standard symplectic case \cite{nek,duist}, the
action-angle coordinates need not be defined globally in $M$, because
of the non-triviality of the fibration $\pi:M\to B$, and are moreover
not unique. But exactly as in that case, they form an atlas in which
any two different sets $(a,\alpha)$ and $(\tilde a,\tilde \alpha)$ of
coordinates are related to each other, in any nonempty connected
component of the intersection of their domains, by a transformation
of the form
\begin{equation}
\label{CambioAA}
  \tilde a=Za+z \,,\qquad \tilde \alpha=Z^{-T}\alpha + \cG(a)
  \ (\mathrm{mod} 2\pi) 
\end{equation}
with $Z\in SL_\pm(n,\interi)$ (the group of unimodular $n\times n$
matrices with integer entries), $z\in\renne$
and an invertible map~$\cG$ (see \cite{FS1}). The only difference from
the symplectic case is that the maps $\cG$ need not satisfy some
conditions imposed by the closedness of $\s$.

A consequence of the semi-global expression \for{eq:sigma} of $\s$ is
that $d\s$ is a basic form. This fact plays a role in
the forthcoming analysis. First of all, since when $n=2$ any basic
2-form is closed, this implies that the case $n=2$ is special:

\begin{proposition}
\label{p:n=2}
If an almost-symplectic manifold $(M,\s)$ with $2$ degrees of
freedom hosts a fully-Lagrangian toric fibration, then $\s$ is
symplectic. 
\end{proposition}

From now on we will thus focus on the case $n\ge3$. Moreover, we
will often regard $d\s$ as defined in $B$.

\vspace{3ex}
\begin{remarks}
\begin{itemize}
\item[(i)] According to \for{eq:sigma}, semi-globally, the
almost-symplectic form is the sum of a symplectic form
($\sum_ida_i\wedge d\a_i$) and of a basic form
($\frac12 \sum_i A_{ij}(a)da_i\wedge da_j$). However, such a
decomposition---and particularly such a symplectic form---might not exist globally. The reason is that, for
the $\sum_ida_i\wedge d\a_i$ to be the local representatives of a 2-form in $M$,
it is necessary that the functions $\cG$ in \for{CambioAA} satisfy
some additional properties (they must be gradients).

\item[(ii)] Reference \cite{FS1} considers a more general situation, which
extends from the symplectic to the almost-symplectic context not only
the notion of complete integrability, but also that of
``noncommutative integrability'' (that allows for isotropic---and not
only Lagrangian---tori). It might have some interest to extend the
analysis of the present article to such more general situation,
because in that case the fibration by tori is typically produced not
by a torus action, but by the action of a non-abelian symmetry group.
\end{itemize}
\end{remarks}

\subsection{Vertical and non-vertical Hamiltonian-like vector fields on
fully-Lagrangian toric fibrations}
\label{ss:3.2}
The existence of a fully-Lagrangian toric fibration $\pi:M\to B$ on an
almost-symplectic manifold $(M,\s)$ forces special properties on the
partially- and fully-Hamiltonian vector fields on $M$. 

First, partially-Hamiltonian vector fields on a generic
almost-symplectic manifold need not conserve the volume induced by the
almost-symplectic form. Instead:

\begin{proposition}
\label{p:volume-2}
If an almost-symplectic manifold $(M,\s)$ hosts a fully-Lagrangian
toric Lagrangian fibration $\pi$, then all locally
partially-Hamiltonian vector fields on $M$ preserve the volume form
$\sigma^n$.
\end{proposition}

\begin{proof} If $i_X\s$ is closed, then $L_{X} \s^n =
n\s^{n-1}\wedge L_X\s = n\s^{n-1}\wedge i_{X}d\s$. From the
expressions in action-angle coordinates of $X$ and $\s$ one sees that
the local representative of $i_Xd\s$ is a sum of terms each of which
contains the wedge product of at least $2$ factors $da_i$
($i=1,\ldots,n$). Hence, $\s^{n-1}\wedge i_Xd\s$ is a sum of terms
each of which contains the wedge product of at least $n+1$ factors
chosen among $da_1,\ldots,da_n$, and therefore vanishes. \end{proof}

Next, all {\it vertical} partially-Hamiltonian vector fields are
automatically fully-Hamiltonian and also integrable (in the sense that
their flow is conjugate to a linear flow on tori).
Indeed, from Lemma \ref{p:basicfunctions-1}, if $X_f^\s$ is
vertical then $f$ is basic and $X_f^\s$ is fully-Hamiltonian;
moreover, the local representative $\cF$ of $f$ in any system
of action-angle coordinates $(a,\a)$ is independent of the angles
and, from \for{eq:HVF},
$$
  (X^\s_f)_{a_i}=0 , \qquad (X^\s_f)_{\a_i}=\der\cF{\a_i}(a)
  \,,\qquad i=1,\dots,n \,.
$$
In addition, in the domain of any action-angle chart, any vertical
partially-Hamiltonian vector field $X^f_\s$ is also Hamiltonian with
respect to the symplectic 2-form $\sum_ida_i\wedge d\a_i$. 

Hence, semi-globally, vertical partially-Hamiltonian vector fields
are standard completely integrable symplectic-Hamiltonian systems.
Globally, however, there might not exist a
symplectic structure on $M$ with such local representatives (this
situation should be understood better).

On the contrary, there may exist both partially-Hamiltonian and
fully-Hamiltonian {\it non-vertical} vector fields. Our aim is to
investigate the structure of the fully-Hamiltonian ones.

\subsection{Fully-Lagrangian toric fibrations as principal
$\tenne$-bundles}
\label{ss:TorusActions}
In the domain $V$ of any of its action-angle charts, any
fully-Lagrangian toric fibration $\pi:M\to B$ is a principal bundle
under the $\tenne$-action on $V$ given by translation of the angles
($\Psi_\beta(a,\a)=(a,\a+\beta))$. This action is fully-Hamiltonian:
given any $\xi\in\renne$, the infinitesimal generator
$Y_\xi=\sum_i\xi_i\partial_{\a_i}$ equals $X_f^\s$ with the basic
function $f=\sum_i\xi_ia_i$. We will consider the case in which this
local action can be extended to a fully-Hamiltonian $\toro^n$-action
on all of $M$. Such a case can be characterized in terms of local
systems of action-angle coordinates:

\begin{definition}
\label{d:GlobalActions} 
A fully-Lagrangian toric fibration $\pi:M\to B$ has {\it global
actions} if it has an action-angle atlas $\cP$ with the following
properties:
\begin{itemize}
\item[i.] All transition functions among the charts of $\cP$ have the
action component equal the identity (all matrices $Z=\unity$ and
vectors $z=0$ in \for{CambioAA}). 
\item[ii.] The map $\hat a^\cP : M \to \renne$
whose value at a point $m\in M$ equals the value at $m$ of the
action coordinates in one (and hence any) chart of $\cP$ whose domain
contains $m$, has connected fibers.
\end{itemize}
The atlas $\cP$ will be called {\it atlas with global actions} and
the map
\beq{a_P}
  \hat a^\cP:M\to \hat a^\cP(M) \subseteq \renne
\eeq
will be said to be the {\it global action map} relative to it.
\end{definition}

Under the hypotheses of this Definition, the global action map $\hat
a^\cP$ as in \for{a_P} is a surjective submersion whose fibers are
connected and diffeomorphic to $\tenne$. By the Ehresman fibration
theorem, $\hat a^\cP$ is a fibration and its fibers, being connected,
coincide with the fibers of~$\pi$. Thus, $\hat a^\cP:M\to \hat
a^\cP(M)$ and $\pi:M\to B$ are different descriptions of the same
fully-Lagrangian toric fibrations, and are related by the fact that
the map
$$
  \hat a_B^\cP:B\to \hat a^\cP(M) 
$$
defined by $\hat a^\cP=\hat a^\cP_B\circ\pi$ is a diffeomorphism.
Considering $\hat a^\cP$ instead of $\pi$ has the advantage that
$\hat a^\cP$ can be viewed as the (equivariant) momentum map of a torus action:

\begin{proposition}
\label{p:converse} Let $(M,\s)$ be an almost-symplectic manifold of
dimension $2n$.
\begin{itemize}
\item[i.] Assume that $\pi:M\to B$ is a fully-Lagrangian toric
fibration with global actions. If $\cP$ is an atlas with
global actions, then $\hat a^\cP:M\to \hat
a^\cP(M)$ is the momentum map of a fully-Hamiltonian free action
$\Psi^{\toro^n}$ of $\toro^n$ on $M$. In any chart of $\cP$,
$\Psi^{\toro^n}$ acts by translation of the angles.

\item[ii.] Consider a fully-Hamiltonian free action $\Psi$ of
$\toro^n$ on $M$ which has a momentum map $J$ with connected fibers.
Then, $J:M\to J(M)\subseteq\renne$ is a fully-Lagrangian toric
fibration with global actions. In the chart of any atlas
with global actions of this fibration, $\Psi$ acts by translation of
the angles.
\end{itemize}
\end{proposition}

\begin{proof} 
(i.) This follows from the remark at the beginning of the section:
for any $\xi\in\renne$, the local representatives of $\sum_i\xi_i
\partial_{\a_i}$ of the infinitesimal generator $Y_\xi$ in the
action-angle charts of $\cP$ are the local representatives of
$X^\s_{J_\xi}$ with $J_\xi=\sum_i\xi_i \hat a^\cP_i$. Since $J_\xi$
is basic, $X^\s_{J_\xi}$ is fully-Hamiltonian. Thus, $\Psi^{\toro^n}$
is a  fully-Hamiltonian action, with momentum map $\hat a^\cP$.

(ii.) Since the action is free, its momentum map $J$ is a submersion
and its fibers are $n$-dimensional submanifolds. Since the group is
abelian, the fibers of $J$ are $\Psi$-invariant. As such, each of them
is a disjoint union of $\Psi$-orbits and, being connected, it is an orbit
of $\Psi$. On the other hand, a free action of a compact group
defines a principal bundle. Thus, $J:M\to J(M)$ is a principal
$\tenne$-bundle; in particular, it is a fibration with compact and
connected fibers.

Choose $\eta,\xi\in\reali^n$. The infinitesimal generator
$Y_\eta=X^\s_{J_\eta}$ is tangent to the fibers of $J$ (because they
are orbits of $\Psi$) and $J_\xi$ is constant on the fibers of $J$.
Therefore, $L_{Y_\eta}J_\xi=0$ and so $\s(Y_\eta,Y_\xi)=
L_{Y_\eta}J_\xi=0$. This proves that the fibers of $J$ are Lagrangian
submanifolds.

Since $J:M\to J(M)$ is a principal $\tenne$-bundle, there are local
trivializations $(J,\phi):U\to \renne\times\tenne$ in which the
$\Psi$-action is given by translation of the angles $\phi$. For each
$i=1,\ldots,n$, the infinitesimal generator
$Y_{e_i}=\partial_{\phi_i}$, where $e_i$ the $i$-th canonical basis
vector of $\mathfrak t^n=\renne$, coincides with $X^\s_{J_i}$. (Here,
$J_i$ is the $i$-th component of the map $J$, and coincides with the
Hamiltonian $J_{e_i}$). 

Consider a function $\tilde f:J(M)\to \reali$ and let $f=\tilde
f\circ\pi$. In a local trivialization $(J,\phi)$, the
partially-Hamiltonian vector field $X_f^\s =
\sum_i\der{\tilde f}{J_i}X^\s_{J_i}$ is a linear combination of
fully-Hamiltonian vector fields and hence is fully-Hamiltonian. Thus, all
basic functions are fully-Hamiltonian. 

This proves that $J:M\to J(M)$ is a fully-Lagrangian
toric fibration. It remains to prove that it has global 
actions. Consider an atlas of $M$ formed by local trivializations $(J,\phi)$ of
the $\tenne$-bundle $J:M\to J(M)$. In any one of them, the local
representative of the almost-symplectic 2-form can be written as
$$
  \sum_{hk}P_{hk}(J,\phi)dJ_h\wedge d\phi_k + 
  \frac12\sum_{hk}Q_{hk}(J,\phi)dJ_h\wedge dJ_k + 
  \frac12\sum_{hk}R_{hk}(J,\phi)d\phi_h\wedge d\phi_k
$$ 
with certain functions $P_{hk}$, $Q_{hk}$, $R_{hk}$. The fact that the
fibers of $J$ are Lagrangian submanifolds implies that all $R_{hk}$
vanish. The fact that $X^\s_{J_i}=\partial_{\phi_i}$ for all $i$
implies $P_{hk}=\delta_{hk}$ for all $h,k$. And the fact that
$L_{\partial_{\phi_i}}\s = L_{X_{J_i}^\s}\s = 0$ for all $i$ implies
that all $Q_{hk}$ are independent of the angles~$\phi$. Hence, the
local representatives of $\s$ have the form \for{eq:sigma}, and these
local trivializations form an action-angle atlas. The local actions
$J_1,\ldots,J_n$ are the restriction to the chart domains of  the
components of the momentum map $J$. From this, and the fact that $J$
has connected fibers, it follows that this atlas has global actions.
\end{proof}

\begin{remark} It follows from this proof that, in the situation of
item ii. of Proposition \ref{p:converse}, any atlas for $M$ which is
formed by local trivializations $(J,\phi)$ of the $\tenne$-bundle
$J:M\to J(M)$ with the property that $\tenne$ acts by translation of
the angles $\a$, is also an atlas with action-angle coordinates for
$M$. We will need this remark later. 
\end{remark}

\section{Symplectization of non-vertical fully-Hamiltonian vector
fields}
\label{s:4}

Constraints on the structure of non-vertical fully-Hamiltonian vector
fields on fully-Lagrangian toric fibrations come from Proposition
\ref{p:upperbound}, according to which, under suitable conditions, a
fully-Hamiltonian function is locally independent of at least $3$
coordinates. Our first goal will be to prove that, under suitable
conditions, these coordinates can be chosen to be angles. Next, we
will investigate the consequences of this fact. 

\subsection{Fourier series and Fourier genericity}
\label{ss:Ipotesi_per_n=4-2}  
Our investigation rests on some assumptions on the non-basic
fully-Hamiltonian functions that we consider. Intuitively, we want
these functions to be ``as less basic as possible'' or, at a
local level, to depend ``as much as possible'' on the angles. 
We will consider two variants of a genericity condition introduced
by Poincar\'e in his study of the non-existence of first integrals in
nearly integrable Hamiltonian systems (\cite{poincare}, Volume 1,
Chapter 5; see also \cite{bfgg}), which refer to the abundance of
harmonics of a function. In order to formulate these conditions we
need a Fourier series on the fibration $\pi:M\to B$, which exists if
$\pi$  is a $\tenne$-principal bundle (see
\cite{FF-hpt,FF-cg,FFS}).

Consider a fully-Lagrangian toric fibration $\pi:M\to B$ with global
actions. Choose an atlas $\cP$ with global actions. Every function
$f\in\cC^\infty(M)$ can be uniquely written as a series
\begin{equation}
\label{FourierSeries}
  f = \sum_{\nu\in\interi^n}f^\nu \,,
\end{equation}
where the functions $f^\nu\in\cC^\infty(M)$, that we will call the
{\it harmonics of $f$ relative to $\cP$}, are defined as follows. The
representative $\cF(a,\a)$ of $f$ in any chart of $\cP$ can be
expanded in the Fourier series
$
  \cF= \sum_{\nu\in\interi^n} \cF^\nu
$,
with ``harmonics''
\beq{harmonics}
  \cF^\nu(a,\a) := \frac{e^{i\nu\cdot\a}}{(2\pi)^n}\int_{\tenne}
  \cF(a,\beta)e^{-i\nu\cdot\beta}d\beta \,.
\eeq
Using \for{CambioAA} with $Z=\unity$ and $z=0$ it is immediate to
verify that, for each $\nu\in\interi^n$, the $\nu$-th harmonics of
the local representatives match in the intersection of the chart
domains, and hence define a function $f^\nu$ on $M$. Clearly, the
series \for{FourierSeries} converges absolutely on any set
$\pi^{-1}(K)$ with $K\subset B$ compact.

The harmonics of a function defined in this way depend on the atlas
only for their labelling with the integer vectors $\nu$. Specifically,
if $\cP$ and $\cP'$ are two atlases with global actions, then in each
connected component of $B$ there exists a matrix $Z\in
SL_\pm(n,\interi)$ such that the harmonics in the two atlases of the
restriction of $f$ to that set are related by
\beq{fourier}
   f^{\nu}_{\cP'}= 
   \det(Z) \, f^{Z\nu}_{\cP} \qquad \forall\nu\in\interi^n \,. 
\eeq

We say that two nonzero integer vectors $\nu,\tilde\nu\in\interi^n$
are {\it parallel} if $\nu=q\tilde\nu$ for some
$q\in\razionali\setminus\{0\}$. Moreover, given a function
$f\in\cC^\infty(M)$, for any $\nu\in\interi^n\setminus\{0\}$ we
define the set
\begin{equation}
\label{B-nu}
  B^{\nu}(f) :=
  \big\{ b\in B\,:\; 
  \exists\ \tilde\nu\not=0 \mathrm{\ parallel\ to\ }
  \nu \ \mathrm{s.t.}\ 
  f^{\tilde\nu}|_{\pi^{-1}(b)}\not=0 \big\} \,. 
\end{equation}
Clearly, $B^\nu(f)=B^{\tilde\nu}(f)$ if $\nu$ and $\tilde\nu$ are
parallel.

\begin{definition}\label{d:PoincGen}
Consider a fully-Lagrangian toric fibration $\pi:M\to B$ with global
actions. A function $f\in\cC^\infty(M)$:
\begin{itemize}
\item[i.] {\it Has property FG1} if, for each
$\nu\in\interi^n\setminus\{0\}$, $B^\nu(f)\subseteq B$ is either empty or open
and dense in $B$.
\item[ii.] {\it Has property FG2} if, for each
$\nu\in\interi^n\setminus\{0\}$, either $B^\nu(f)\subseteq B$ is empty or its
complement is countable.
\end{itemize}
\end{definition}

Note that only the labelling of the sets $B^\nu$ depends on the choice of the
atlas, and the fact that a function has either properties FG1 or
FG2 is independent of this choice. 

\begin{remarks} 
\begin{itemize}
\item[(i)] The fact that properties FG1 and FG2 allow for the
possibility that all harmonics parallel to a given one are zero is
certainly not a genericity condition in the class of all smooth
functions. Allowing for it is however unavoidable
because---essentially as a consequence of
Proposition~\ref{p:upperbound}---the representatives of
fully-Hamiltonian functions cannot depend on all angles. Thus, these
properties should be thought of as genericity conditions in the class
of fully-Hamiltonian functions. 

\item[(ii)] The construction of the Fourier series does not need that the
fibration $\pi$ has global action, but only that there is an
action-angle atlas whose all transition functions have matrices
$Z=\unity$ (one such fibration could be said to have ``no monodromy''
by analogy with the symplectic case \cite{nek,duist}).
\end{itemize}
\end{remarks}

\subsection{Symplectization of non-vertical fully-Hamiltonian
vector fields}
\label{ss:results} 
We state now our results on the existence and structure of
non-vertical fully-Hamiltonian vector fields on fully-Lagrangian
toric fibration. We will assume that $M$ (and hence $B$) is connected;
if not, all results apply to each component separately.

As already remarked, cases with $n=1$ and $n=2$ are
trivial (see Proposition \ref{p:n=2}). The cases $n=3$ and $n=4$ are
special, and we can give results specific to them. 

First, generically, under mild genericity conditions on either $\s$
or the Hamiltonian, when $n=3$ there are no such vector fields: all
fully-Hamiltonian vector fields are vertical and hence,
semi-globally, Hamiltonian and completely integrable with respect to
the semi-global symplectic structure $\sum_i da_i\wedge d\alpha_i$.

\begin{theorem}
\label{thm1}
Let $(M,\s)$ be a connected almost-symplectic manifold with $n=3$
degrees of freedom that hosts a fully-Lagrangian toric fibration
$\pi:M\to B$. 
\begin{itemize}
\item[i.] If $d\s$ is nonzero in an open and dense subset of
$B$, then any fully-Hamiltonian vector field is vertical.
\item[ii.] Assume that $\pi$ has global actions and that $d\s\not=0$.
If $f\in\cH_F^\s(M)$ has property FG1, then $X^\s_f$ is vertical.
\end{itemize}
\end{theorem}

\noindent
Instead, when $n=4$, non-vertical fully-Hamiltonian vector fields may
exist. However, under the genericity condition FG1, they
are Hamiltonian and completely integrable with respect to {\it some}
semi-global symplectic structure:

\begin{theorem}
\label{thm3}
Let $(M,\s)$ be a connected
almost-symplectic manifold with $n=4$ degrees of freedom and
$d\s\not=0$ which hosts a fully-Lagrangian toric fibration $\pi:M\to
B$. If $f\in\cH_F^\s(M)$ has property FG1, then
semi-globally, in a neighbourhood $V$ of each fiber of $\pi$,
$X^\s_f$ is Hamiltonian, and completely integrable, with respect to a
symplectic structure $\s^s_V$ on $V$. 
\end{theorem}

Our next result is that, under the stronger genericity condition FG2,
when $n\ge5$ all non-vertical fully-Hamiltonian vector fields reduce,
under almost-symplectic torus actions, to families of standard
Hamiltonian vector fields on reduced {\it symplectic} manifolds. 
This is true also when $n=3,4$, and we thus state this result here for
all $n\ge3$; the connection with Theorems \ref{thm1} and \ref{thm3}
will be discussed in Section~\ref{s:esempi}. Remember that we make
the convention that a zero-dimensional almost-symplectic manifold is
a symplectic manifold. 
 
\begin{theorem}
\label{thm2}
Let $(M,\s)$ be a connected
almost-symplectic manifold with $n\ge3$ degrees of freedom and
$d\s\not=0$ which hosts a fully-Lagrangian toric fibration $\pi:M\to
B$ with, if $n\le5$, global actions. 
Consider a function $f\in\cH_F^\s(M)$ which
has property FG1 and, if $n\ge 6$, property FG2. 

Then, globally in $M$
if $n\le5$, and in a neighbourhood $V$ of each fiber of
$\pi$ (namely, semi-globally) if $n\ge6$, there exist a 
$3\le k\le n$ and a fully-Hamiltonian action $\Psi^{\toro^k}$ of
$\toro^k$ on $M$ such that:
\begin{itemize}
\item[i.] $f$ is $\Psi^{\toro^k}$-invariant.
\item[ii.] On each almost-symplectic reduced manifold $(V_c,\s_c)$
there is a symplectic structure $\s_c^s$ such that
$X^{\s_c}_{f_c}=X^{\s_c^s}_{f_c}$.
\end{itemize}
(Here, $(V_c,\s_c)$ is the reduced space, in the sense of Proposition
\ref{p:vaisman}, of the almost-symplectic manifold $(V,\s|_V)$ under
the action $\Psi^{\toro^k}$).
\end{theorem}

\noindent
We will discuss some aspects of the {\it dynamics} of these vector
fields, and give some examples of their existence, in Section 6. We
anticipate that the following situations are possible, depending on
the number $r:=n-k$ of degrees of freedom of the reduced symplectic
spaces:
\begin{itemize}
\item[1.] If $r=0$, then all fully-Hamiltonian vector fields are
vertical (as in the case $n=3$ of Theorem \ref{thm1}).
\item[2.] If $r=1$, then all fully-Hamiltonian vector fields are
completely integrable and Hamiltonian with respect to some
semi-global symplectic forms (as in the case $n=4$ of Theorem
\ref{thm3}).
\item[3.] If $r\ge2$, then there can be a variety of situations; 
in particular, non-vertical fully-Hamiltonian vector fields need not
be integrable and may have chaotic dynamics.
\end{itemize}

\begin{remarks} 
\begin{itemize}
\item[(i.)] The $\toro^k$-action in Theorem \ref{thm2} is the
action of a subgroup of the torus $\tenne$ that generates the
fully-Hamiltonian toric fibration. Correspondingly, there exists an
action-angle atlas with global actions with the property that the
local representative of $f$ in each chart are independent of the first $k$
angles. 

\item[(ii.)] Since Theorem \ref{thm1} does not involve any group action, for
its validity it would be sufficient to require the existence of an
atlas with all changes of actions equal the identity (that is, ``no
monodromy''), without requiring the connectedness of the fiber of the
map~$\hat a^\cP$.
\end{itemize}
\end{remarks}

\section{Proof of Theorems \ref{thm1}, \ref{thm3} and \ref{thm2} }
\label{s:proofs}

At the basis of the three Theorems stands the estimate of Proposition
\ref{p:upperbound}. Item i. of Theorem \ref{thm1} follows more or
less directly from it (and from a simple consequence of it, Lemma
\ref{l:kerC} below). The proof of the other statements requires
additional considerations. The key fact will be to prove that, under
the hypotheses there made, there is an atlas with global actions in
which all representatives of fully-Hamiltonian functions are
independent of the same 3 angles (Lemma \ref{l:scelta-angoli}).
Translation of these angles gives the torus action. This directly
implies item ii. of Theorem \ref{thm1} and Theorem \ref{thm3}. The
proof of Theorem \ref{thm2} uses an iteration of this procedure.

\subsection{The local representatives of fully-Hamiltonian functions}
\label{ss:Spectra}
Consider a fully-Lagrangian toric fibration $\pi:M\to B$ of an almost
symplectic manifold $(M,\s)$ (for now, we need not require that $M$
is connected and that $\pi$ has global actions) and an action-angle
chart $(a,\a):V\to\cA\times\tenne$. The representative $\saa$ of $\s$
in this chart is as in \for{eq:sigma}. Hence 
\begin{equation}
\label{eq:dsigma}
  d\saa \;=\; 3\sum_{ijk} C_{ijk} da_k\otimes da_i\otimes da_j 
\end{equation}
with 
\begin{equation}\label{eq:C}
  C_{ijk}(a) := \der{A_{ij}}{a_k}(a) + \der{A_{ki}}{a_j}(a) 
  + \der{A_{jk}}{a_i}(a) 
  \,,\qquad i,j,k=1,\ldots,n 
\end{equation}
(unless otherwise said, all summation indexes run from 1 to $n$).
The condition that $\s$ is not closed in $V$ is precisely that 
the skew-symmetric 3-tensor field $C$ with components $C_{ijk}$ does
not vanish. (This tensor is the almost-symplectic analogue of the
``Jacobiator'' of almost-Poisson manifolds \cite{cannas, koon,
balseiro1}).

A function $\cF$ on $\cA\times\tenne$ is the local representative of
a fully-Hamiltonian function on $M$ if and only if
$i_{X_\cF^{\saa}}d\saa = -\sum_{ijk}C_{ijk}\der \cF{\a_i} da_j\otimes
da_k$ vanishes, namely if and only if
\begin{equation}
\label{Cderfa=0}
  C(a)\der \cF\a(a,\a) =  0 \qquad \forall (a,\a) \in\cA\times \tenne
  \,.
\end{equation}
Here, $C(a)u=0$ means of course $\sum_kC_{ijk}(a)u_k=0$ for all
$i,j=1,\ldots,n$.

In agreement with the fact that $d\s$ is a basic 3-form, we will
regard the tensor field $C$ as defined on $\cA$. At each point
$a\in\cA$,
$$
   \ker C(a) \,:=\, 
   \big\{u\in\renne \,:\; C(a)u=0 \big\} 
$$
is thus a linear subspace of $\renne$.

We define the {\it spectrum} $\Sp{\cF,a}$ of a function
$\cF:\cA\times\tenne\to\reali$ at a point $a\in \cA$ as
\begin{equation}
\label{spectrum2}
  \Sp{\cF,a} :=
  \big\{ \nu\in\interi^n\setminus\{0\} 
  \,:\; \widehat\cF(a) \not=0 \big\} 
\end{equation}
where, for each $\nu$, the ``Fourier coefficient'' $\widehat\cF^\nu$
is defined by $\cF^\nu(a,\a)=:\widehat\cF^\nu(a)e^{i\nu\cdot\a}$ with
$\cF^\nu$ the harmonics \for{harmonics}. Note that we have removed the
zero-harmonics $\cF^0$ from the spectrum.

We denote by $\big\langle U\big\rangle$ the linear subspace of
$\renne$ generated by a subset $U\subset\renne$, namely, the smallest
subspace of $\renne$ that contains $U$.

\begin{lemma}
\label{l:kerC}
Consider an action-angle chart $(a,\a):V\to\cA\times\tenne$ of 
a fully-Lagrangian fibration on an almost-symplectic manifold with
$n\ge3$ degrees of freedom.
\begin{itemize}
\item[i.] A function $\cF$ on $\cA\times\tenne$ is the local
representative of a fully-Hamiltonian function 
if and only if
\begin{equation}
\label{spectrum-2}
  \big\langle\Sp{\cF,a}\big\rangle \;\subseteq\; \ker C(a) 
  \qquad \forall a\in\cA \,.
\end{equation}
\item[ii.] At a point $a\in\cA$ at which $C(a)\not=0$, $\ker C(a)$
is a subspace of $\renne$ of dimension $\le n-3$.
\end{itemize}
\end{lemma}

\begin{proof} (i.) Expanding $\cF$ in Fourier series, the condition
of full-Hamiltonianity \for{Cderfa=0} becomes
$$
 \sum_{\nu\in{\interi^n}} \sum_k C_{ijk}(a)\nu_k \cF^\nu(a,\a) =0 
 \qquad \forall i,j,a,\a \,,
$$
namely
$$
 \sum_k C_{ijk}(a)\nu_k \widehat\cF^\nu(a) = 0 
 \qquad \forall i,j,a,\nu \,.
$$
This is equivalent to the fact that
$C(a)\nu=0$ for all $\nu$ and $a$ such that $\widehat\cF^\nu(a)
\not=0$, namely, see \for{spectrum2}, 
$$
   \Sp{\cF,a}\subset \ker C(a) \quad \forall a \,.
$$
Since $\ker C(a)$ is a linear subspace of $\renne$, this in turn
implies $\big\langle \Sp{\cF,a}\big\rangle \subseteq \ker C(a)$.

(ii.) See the proof of item ii. of  Proposition~\ref{p:upperbound}.
\end{proof}

\subsection{Proof of item {i.} of Theorem \ref{thm1}}
\label{ss:ProofThm1.i}
Let $\cF$ be the representative of a fully-Hamiltonian function $f$ in
an action-angle chart $(a,\a):V\to\cA\times\tenne$. Under the
hypotheses of item i. of Theorem \ref{thm1}, the tensor $C$ is different from
zero in an open and dense subset $\cA_C$ of $\cA$. Since $n=3$, by
Lemma \ref{l:kerC} the kernel of $C$ is zero-dimensional at all
points of $\cA_C$. Hence condition \for{Cderfa=0} implies $\der
\cF\alpha(a,\alpha)=0$ at all points $(a,\a)\in \cA_C\times\toro^n$.
By continuity, $\der \cF\alpha=0$ in all of $\cA\times\toro^n$ and
$\cF$ is independent of the angles. By the arbitrariness of the chart,
this proves that $f$ is basic and, by Lemma
\ref{p:basicfunctions-1}, $X_f^\s$ is vertical.

\subsection{The spectrum of fully-Hamiltonian functions}
For a fully-Lagrangian toric fibration $\pi:M\to B$ with global
actions it is possible to extend the definition of the spectrum to
functions on $M$ (the connectedness of $M$ is not needed for this, but
will be needed in the next Lemma).

Specifically, choose an action-angle atlas $\cP$ with global actions,
and use the Fourier series on $\cC^\infty(M)$ relative to it. Define
the {\it spectrum} (relative to $\cP$) of a function
$f\in\cC^\infty(M)$ at a point $b\in B$ as the set 
\begin{equation}
\label{spectrum}
  \Sp{f,b} \,:=\,
  \big\{ \nu\in\interi^n\setminus\{0\} 
  \,:\; f^\nu\big|_{\pi^{-1}(b)} \not=0 \big\} \,.
\end{equation}
The spectrum defined in this way depends on the atlas $\cP$, but the
dependence is rather trivial: if $\cP'$ is another atlas with global
actions, then, by \for{fourier}, in each connected component of $B$
there exists a unimodular $n\times n$ integer matrix $Z$ which, for
any function $f$ and point $b$, relates the spectra $\Sp{f,b}$
relative to $\cP$ and $\mathrm{Sp}'(f,b)$ relative to $\cP'$,
$\mathrm{Sp'}(f,b)=Z\mathrm{Sp}(f,b)$.

The proofs of items~ii.\@ of Theorems \ref{thm1} and \ref{thm2} are
based on the choice of a suitable action-angle atlas, as in the
following Lemma:

\begin{lemma}
\label{l:scelta-angoli}
Consider a connected almost-symplectic manifold $(M,\s)$ with $n\ge3$
degrees of freedom and $d\s\not=0$ which hosts a fully-Lagrangian
toric fibration $\pi:M\to B$ with global actions. For any function
$f\in\cH_F^\s(M)$ with property FG1 there exists an atlas $\cP^f$ of
$M$ with global actions in which all local representatives of $f$ are
independent of the last $3$ angles.
\end{lemma}

\begin{proof}
Choose an atlas $\cP$ with global actions and use the Fourier series
\for{FourierSeries} and the spectra \for{spectrum} relative to it.
Define
$$
   \Sp{f}:=\bigcup_{b\in B} \Sp{f,b} \subseteq \interi^n \,.
$$
Since $f$ has property FG1, for each $\nu\in\Sp f$ the set $B^\nu(f)$
as in \for{B-nu} is an open and dense subset of $ B$. Since $ B$
is a differentiable manifold, and hence a Baire space, the countable
intersection 
$$
  B^*(f):= \bigcap_{\nu\in\Sp f} B^\nu(f) 
$$
is a dense subset of $B$ \cite{kelley}. 

Given that $d\s$ is continuous, the subset $\cO$ of $B$
where $d\s\not=0$ is open (recall that since $d\s$ is a basic form we
can regard it as defined in $B$). Since $\s$ is not closed, $\cO$ is
nonempty and $\cO\cap B^*(f)$ is nonempty as well.

Choose a point $b^*\in\cO\cap B^*(f)$. By the property
FG1 of $f$, for any $\nu\in\Sp f$ there is a $\tilde\nu\not=0$
parallel to $\nu$ which belongs to $\Sp{f,b^*}$. (Recall that
$0\notin\Sp f$). Therefore, 
$\Sp f \subset \big\langle\Sp{f,b^*}\big\rangle \subseteq\renne$ and so 
$$
  \big\langle\Sp f\big\rangle =
  \big\langle\Sp{f,b^*}\big\rangle \,.
$$
Since the Fourier series we are using is relative to $\cP$,
$\Sp{f,b^*}$ equals the spectrum,  at the corresponding point
$a^*\in\cA$, of the representative of $f$ in any chart of $\cP$ that
contains $b^*$. Since $d\s(b^*)\not=0$, it follows from Lemma
\ref{l:kerC} that $\big\langle\Sp{f,b^*}\big\rangle$ is a linear
subspace of $\ker C(a^*)$ and hence has dimension $\le n-3$.
Therefore, 
$$
 \cL := \big\langle\Sp{f}\big\rangle \cap\interi^ n
$$
is a sublattice of $\interi^n$ of rank $r \le n-3$. (A sublattice of
$\interi^n$ of rank $r$ is the set of all the linear combinations
with integer coefficients of $r$ vectors $u_1,\ldots,u_r\in
\interi^n$, called a basis of the sublattice).

If $r=0$, then $\Sp{f}=\{0\}$ and all representatives of $f$ are
independent of all angles. Assume $r\ge1$.

The lattice $\cL$ is not just a generic sublattice of $\interi^n$,
but it is the intersection of $\interi^n$ with a linear subspace of
$\renne$. This guarantees that any basis of $\cL$ can be completed to
a basis of $\interi^n$. [The Elementary Divisors Theorem
(\cite{lang}, Theorem 7.8) guarantees that for any finitely generated
submodule $\not=\{0\}$ (e.g., a lattice of positive rank) of a free
abelian module over a principal ideal domain (e.g., $\interi^n$)
there exists a basis $\{u_1,\ldots,u_n\}$ of the latter, an integer
$1\le r\le n$ and integers $d_1,\ldots, d_r$ such that
$\{d_1u_1,\ldots,d_ru_r\}$ is a basis of the former. In our case, all
$d_i=1$ because the lattice $\cL$ is the intersection of $\interi^n$
with a subspace of $\renne$.]

Choose a basis $\{u_1,\ldots,u_r\}$ of $\cL$ and complete it to
a basis $\{u_1,\ldots,u_r,u_{r+1},\ldots,u_n\}$ of $\interi^n$.
Then, there exists a unimodular integer matrix $Z$ such that
$Zu_i=e_i$, the $i$-th unit vector, for all $i=1,\ldots,n$. 
The change of coordinates
$$
  \cC: \cA\times\tenne \to \tilde\cA\times \tenne \,,
  \qquad (a,\a)\mapsto (\tilde a,\tilde \a) = (Z^Ta, Z^{-1}\a) 
$$
in any chart $(a,\a)\in\cA\times\tenne$ of the considered atlas
produces a new atlas with action-angle coordinates $\cP^f$ for
$M$. The local representatives $\tilde \cF$ of $f$ in the charts of
$\cP^f$ have nonzero Fourier components only for $\nu$ in the
subspace generated by $e_1,\ldots,e_r$. Hence,
$\der{\tilde\cF}{\tilde\a_i}=0$ for all $i=r+1,\ldots,n$ and
$\tilde\cF$ is independent of the last $k=n-r\ge 3$ angles. 

Just as the original atlas, the new atlas $\cP^f$ has all
transition functions between the action coordinates given by the
identity. Given the linearity of the transition functions between the
actions of $\cP$ and $\cP^f$, if $\cP$ has global 
actions then so does $\cP^f$. 

We also note, for future use, that if the atlas $\cP$ consists
of a single chart, then also the atlas $\cP^f$ built in this proof
consists of a single chart. \end{proof}

\subsection{Proof of item ii.\@ of Theorem \ref{thm1}}
\label{ss:ProofThm1.ii}
By Lemma \ref{l:scelta-angoli}, when $n=3$ all local representatives
in a certain action-angle atlas of a function $f\in\cH_F^\s(M)$ which
has property FG1 are independent of the angles. Hence $f$ is basic
and $X_f^\s$ is vertical.

\subsection{Proof of Theorem \ref{thm3}}
\label{ss:ProofThm3}
Choose any one of the action-angle charts provided by Lemma
\ref{l:scelta-angoli}, in which the local representative $\cF$ of $f$
is independent of the last three angles. Denote these coordinates
$
  (a,\a)=(I,J,\phi,\psi) \in
  \reali\times\reali^{3} \times 
  \toro^{1}\times\toro^{3},
$
so that $\cF$ is independent of the angles $\psi$. Restrict the
domain of this chart to an open set $V$ (which is the one appearing in
the statement of the theorem) such that $(I,J,\phi,\psi)\in\mathcal
I\times\mathcal{J}\times\toro^1\times\toro^3$ with an
interval $\mathcal{I}\ni I$ and an open set
$\mathcal{J}\ni J$. 

The representative of $\s$ in this chart has the form
$$
  \saa=dI\wedge d\phi + \sum_idJ_i\wedge d\psi_i
  + \sum_iA_{IJ_i}(I,J)dI\wedge dJ_i
  + \frac12 \sum_{i,j}A_{J_iJ_j}(I,J)dJ_i\wedge dJ_j 
$$
where all sum runs from $1$ to $3$. The vector field $X_\cF^{\saa}$,
seen as an ODE, is given by
\beq{queste}
  \dot I = -\der\cF\phi\,,\qquad 
  \dot J_i = 0 \,,\qquad 
  \dot\phi = \der\cF I\,,\qquad
  \dot\psi_i = \der\cF {J_i}
  + A_{IJ_i}\der\cF\phi
  \,,\qquad i=1,2,3 \,. 
\eeq
Choose an $I_0\in \mathcal{I} \subseteq\reali$ and define, for each
$i$, $\cG_i(I,J) := \int_{I_0}^I A_{I,J_i}(x,J)dx$.
Along the flow of $X_f^\s$, $\frac d{dt}\cG_i=A_{IJ_i}\dot I =
-A_{IJ_i}\der\cF\phi$ and the change of coordinates
\beq{quelle}
 (I,J,\phi,\psi) \mapsto (I,J,\phi,\chi=\psi+\cG(I,J)) \,,
\eeq
with $\cG=(\cG_1,\cG_2,\cG_3)$, conjugates equations \for{queste}
to
$$
  \dot I = -\der\cF\phi\,,\qquad 
  \dot J_i = 0 \,,\qquad 
  \dot\phi = \der\cF I\,,\qquad
  \dot{\chi}_i = \der\cF {J_i} 
  \,,\qquad i=1,2,3 \,,
$$
which are Hamilton equations for the symplectic 2-form
\beq{sigma-3}
  \tilde\s^\mathrm{aa} = dI\wedge d\phi + \sum_idJ_i\wedge d\chi_i \,.
\eeq

The functions $\cF$, $J_1$, $J_2$, $J_3$ are first integrals of these
equations and, if $\der\cF\phi=0$, so is also $I$. Since $f$ has
property FG1, $\der\cF\phi$ either is everywhere zero or is different
from zero in an open and dense subset of the chart domain. In the
first case $I,J_1,J_2,J_3$ are first integrals which are in
involution with respect to $\tilde\s^\mathrm{aa}$ and are everywhere
functionally independent in the chart domain (and
$X_f^\s$ is thus vertical in that domain). In the second case,
$\cF,J_1,J_2,J_3$ are first integrals which are in involution with
respect to $\tilde\s^\mathrm{aa}$ and are functionally independent in
an open and dense subset of the chart domain. 

\begin{remark} This proof uses only the fact that $n-k=1$, not that
$n=4$. Thus, this argument could be used to prove statement 2. after
Theorem \ref{thm2}; see also Section \ref{ss:esempi.2}.
\end{remark}

\subsection{Proof of Theorem \ref{thm2}}
\label{ss:ProofThm2}
Lemma \ref{l:scelta-angoli} can be reformulated by saying that, in
the stated hypotheses, any fully-Hamiltonian function on $M$ with
property FG1 is invariant under a fully-Hamiltonian action of a
three-dimensional torus, which corresponds to translation of the last
three angles in each chart of the atlas $\cP^f$. The proof of Theorem
\ref{thm2} is based on the fact that, for functions with property FG2,
this property is semi-globally inherited by the reduced systems---as long as they
are not symplectic. 

Preliminarily, we state a Lemma, which we will then use with
$k=3$. We make the convention that a zero-dimensional 
manifold is (trivially) almost-symplectic and carries a
(trivial) fully-Lagrangian toric fibration.

\begin{lemma}
\label{l:quoziente} Let $(M,\s)$ be a connected almost-symplectic manifold with
$n\ge3$ degrees of freedom and $d\s\not=0$ which hosts a
fully-Lagrangian toric fibration $\pi:M\to B$ with global actions.
Consider an action-angle atlas $\cP$ of $M$ with global actions and,
for some $1\le k\le n$, the action $\Psi^{\toro^k}$ of $\toro^k$ on $M$
by translation of the last $k$ angles in all charts of $\cP$.
\begin{itemize}
\item[i.] $\Psi^{\toro^k}$ is fully-Hamiltonian and, denoting
$J_{\toro^k}$ its momentum map, for each ${c}\in J_{\toro^k}(M)$ the
reduced manifold $(M_{c}=J_{\toro^k}^{-1}(c)/\toro^k,\s_{c})$ has
dimension $2(n-k)$ and carries a fully-Lagrangian
toric fibration $\pi_c:M_c\to B_c$ with global actions.
\item[ii.] If $f\in \cH_F^\s(M)$ is $\Psi^{\toro^k}$-invariant and has
property FG2, then, for each $c\in J_{\toro^k}(M)$,
the reduced function $f_c \in \cH_F^{\s_c}(M_c)$ has property
FG2 as well.
\end{itemize}
\end{lemma}

\begin{proof} For $k=n$ there is nothing to prove, so we assume
$k<n$.

(i.) 
This is all obvious except, perhaps, the existence of the
fully-Lagrangian toric fibration with global actions of $(M_c,\s_c)$,
whose proof is anyway elementary. We detail it here also because in
this way we introduce some notation that we will use again below.

Denote $p_c:J_{\toro^k}^{-1}(c)\to M_c$ the canonical projection,
$\iota_c:J_{\toro^k}^{-1}(c)\hookrightarrow M$ the immersion and $f_c$ the
reduced function induced by $f$ on $M_c$, which satisfies $f_c\circ
p_c=f\circ\iota_c$ (see Section \ref{ss:ASdiffeos}).

The atlas $\cP$ of $M$ induces a ``quotient action-angle atlas''
$\cP_c$ of $M_c$ as follows. Consider a chart
$$
  \cC = 
  (I,J,\phi,\psi) : U \to 
  \reali^{n-k}\times\reali^k \times
  \toro^{n-k}\times\toro^k 
$$
of $\cP$ whose domain $U$ has nonempty intersection with
$J_{\toro^k}^{-1}(c)$. In this chart, the representative of $\s$ has
expression \for{eq:sigma} and $J:U\to\reali^k$ is the
representative of the momentum map $J_{\toro^k}$. Thus
$\widehat U_c:=J_{\toro^k}^{-1}(c)\cap U$ is the subset where $J=c$,
$$
  \widehat\cC_c=(I,\phi,\psi) : \widehat U_c \to 
  \reali^{n-k}\times \toro^{n-k}\times\toro^k
$$
is a chart of $J_{\toro^k}^{-1}(c)$,
$$
  \cC_c = (I,\phi) : U_c \to \reali^{n-k} \times \toro^{n-k} \,,
$$
with $U_c:=p_c(\widehat U_c)$, is a chart of $M_c$, the
representative of $p_c$ relative to the two charts $\widehat\cC_c$
and $\cC_c$ is $(I,\phi,\psi)\mapsto (I,\phi)$, and the
representative of $\s_c$ relative to $\cC_c$ is, from \for{eq:sigma}
and Proposition \ref{p:vaisman},
$$
  \sum_{i=1}^{n-k}dI_i\wedge d\phi_i + 
  \frac12 \sum_{i,j=1}^{n-k}A_{ij}(I,c)dI_i\wedge dI_j \,.
$$ 
Collecting together all these charts produces an atlas $\cP_c$ of
$M_c$ whose transition functions are inherited from those
\for{CambioAA} of $\cP$ and have the form
\beq{CambioAA-2}
  \tilde I=I \,,\qquad \tilde\phi=\phi+\cG_\phi(I,c)
\eeq
with $\cG_I$ the first $n-k$ components of $\cG$. 

From \for{CambioAA-2} it follows that there is a well defined action
of $\toro^{n-k}$ on $M_c$ which, in every chart of
$\cP_c$, acts by translation of the angles. This action is
fully-Hamiltonian and, by item ii. of Proposition \ref{p:converse},
defines a fully-Lagrangian toric fibration $\pi_c:M_c\to B_c$ with
global actions (whose fibers are the sets $I=\const$). $\cP_c$ is an
action-angle atlas for $\pi_c:M_c \to
B_c$ (see the Remark at the end of Section \ref{ss:TorusActions}). 

(ii.) That $f_c$ is fully-Hamiltonian follows from Proposition
\ref{p:vaisman}. The atlas $\cP$ induces a Fourier series on
$\pi:M\to B$ and the atlas $\cP_c$ introduced in the proof of item i.
induces a Fourier series on $\pi_c:M_c\to B_c$. By the invariance of
$f$, the representative $\cF$ of $f$ in any chart of $\cP$ is
independent of the last $k$ angles and hence its restriction to $J=c$
equals the representative $\cF_c$ of $f_c$ in the corresponding chart
of $\cP_c$, that is
\beq{cF}
  \cF(I,c,\phi,\psi) = \cF_c(I,\phi) \,.
\eeq
This implies that if a harmonic $f^\nu$ of $f$ is nonzero then
necessarily $\nu=(\mu,0)\in\interi^{n-k}\times\interi^k$ and
$$
   f_c^{\mu}\circ p_c = f^\nu\circ \imm_c \,.
$$
Fix $\mu\in\interi^{n-k}$ and assume that
$$
  B^{\mu}_c(f_c) =
  \big\{ b_c\in B_c \,:\; 
  \exists\ \tilde\mu\in\interi^{n-k}\setminus\{0\} \mathrm{\ parallel\ to\ }
  \mu \ \mathrm{s.t.}\ 
  f_c^{\tilde\mu}\big|_{\pi_c^{-1}(b_c)}\not=0 \big\} 
$$
is not empty (see \for{B-nu}). Then, there exist $b_c\in B_c$,
$\tilde\mu$ parallel to $\mu$ and $m_c\in \pi_c^{-1}(b_c)$ such that
$f_c^{\tilde\mu}(m_c)\not=0$. Now, $m_c=p_c(m)$ for some $m\in
J_{\toro^k}^{-1}(c)$, and so $f^{(\tilde\mu,0)}(m) =
f_c^{\tilde\mu}(m_c)\not=0$. This implies that, if $B_c^{\mu}(f_c)$ is
not empty then $B^{(\mu,0)}(f)$ is not
empty as well, and thus, since $f$ has property FG2, the complement of
$B^{(\mu,0)}(f)$ is
countable. Consequently, also the complement of
$B^{\mu}_c(f_c)$ is countable and, by
the arbitrariness of $\mu$, $f_c$ has property FG2. \end{proof}

\begin{remark} 
The base $B_c$ of the fibration $\pi_c:M_c\to B_c$ of Lemma
\ref{l:quoziente}  equals $J_{\toro^{n-k}} (J_{\toro^k}^{-1}(c))$,
where $J_{\toro^{n-k}}$ is the momentum map of the action of
$\toro^{n-k}$ on $M$ by translation of the {\it first} $n-k$ angles
in all charts of $\cP$, and also $J_{c,\toro^{n-k}}(M_c)$, where
$J_{c,\toro^{n-k}}$ is the momentum map of the action of
$\toro^{n-k}$ on $M_c$ by translation of the angles in all charts of
$\cP_c$ (the equality of these two expressions can be proved observing
that the two actions $\Psi^{\toro^k}$ and $\Psi^{\toro^{n-k}}$
commute and using elementary facts from an ``almost-symplectic
reduction in stages'' that mimics the symplectic reduction in stages
for commuting actions described for instance in Chapter 4 of
\cite{stages}).
\end{remark}

{\it Proof of Theorem \ref{thm2}.}
By Lemma \ref{l:scelta-angoli}, there is an action-angle atlas $\cP^f$ of $M$
with global actions such that the local representatives of
$f$ in all its charts are independent of the last 3 angles. We may
thus apply Lemma \ref{l:quoziente} with $k=3$. 

This is enough to prove Theorem \ref{thm2} if $n=3,4,5$, because in
that case all $\Psi^{\toro^3}$-reduced spaces of Lemma
\ref{l:quoziente} have either 0, 1 or 2 degrees of freedom
and are symplectic.

Consider now the case $n\ge6$. If all reduced systems are symplectic,
then the proof is concluded, again with $k=3$. 

Assume, thus, that there exist $c\in J_{\toro^3}(M)$ and $m_c\in
M_c$ such that $d\s_c(m_c)\not=0$. By Proposition \ref{p:vaisman} 
$f_c$ is fully-Hamiltonian and by Lemma \ref{l:quoziente} it
has property FG2. However, as noticed in a Remark at the end of
Section \ref{s:2}, the reduced space $M_c$ might be not connected. In
such a case, Lemma \ref{l:scelta-angoli} can be applied only to the
component $M^\circ_{c}$ of $m_c$. Accordingly, the action-angle atlas
in which the representatives of $f_c$ are independent of the last 3
(of the $n-3$, now) angles and the result about the existence of the
torus action would be valid only in a subset of $M$. Iterating this
procedure might lead again to non-connected reduced spaces and to a
proliferation of further subdivisions of $M$.
So, in general, this procedure leads not to a global torus action on
$M$, but to a family of semi-global actions. We thus prefer to start
from the beginning from an open set in $M$ with the properties that
all the reduced spaces are connected. One such open set is the domain
of an action-angle chart with the property that the actions take
values in a convex subset of $\renne$ (see below). Clearly, one such
chart exists in a neighbourhood of any point of $M$. 

Consider, thus, an action-angle chart $\cC=(a,\a):V\to
\cA\times\tenne$, with $V\subseteq M$, such that $\cA\subseteq
\renne$ is convex and the representative $\cF$ of $f$ is independent
of the last three angles. Write
$a=(I,J)\in\reali^{n-3}\times\reali^3$ and 
$\a=(\phi,\psi)\in\toro^{n-3}\times\toro^3$. The momentum map of the
$\toro^3$-action on $V$ by translations of the angles $\psi$  is the
actions $J:V\to\reali^3$. We now make the following two steps.

{\it Step 1.} Fix a value $c_1$ of $J$. The reduced space
$V_{c_1}:=J^{-1}(c_1)/\toro^3$ has the induced action-angle chart
$\cC_{c_1}=(I,\phi):V_{c_1}\to\cA_{c_1}\times\toro^{n-3}$, which is defined by
the fact that the representative of the projection $p_{c_1}:V\to
V_{c_1}$ relative to the two charts $\cC$ and $\cC_{c_1}$ is
$(I,\phi,\psi)\mapsto(I,\phi)$ (see the proof of Lemma
\ref{l:scelta-angoli}). Since $\cA$ is convex, hence connected,
so is $\cA_{c_1}=\{I\in\renne: (I,c_1)\in\cA\}$. 

By Lemma \ref{l:quoziente} $f_{c_1}$ has property FG2. Since
$V_{c_1}$ (being diffeomorphic to $\cA_{c_1}\times\toro^{n-3}$) is
connected, Lemma \ref{l:scelta-angoli} applies to $f_{c_1}$. More
precisely, as noticed at the end of the proof of that Lemma,
since $V_{c_1}$ is covered by a single action-angle chart
there is a
matrix $\widehat Z\in SL_\pm(n-3,\interi)$ such that the local
representative $\hat\cF_{c_1}$ of $f_{c_1}$ in the action-angle
coordinates
$$
  (\hat Z^TI,\hat Z\phi) : 
  V_{c_1} \to (\hat Z^T\cA_{c_1})\times\toro^{n-3}
$$
on $V_{c_1}$ is independent of the last three angles.

{\it Step 2.} Let $Z=\mathrm{diag}(\widehat Z,\unity_3)$, where
$\unity_3$ is the $3\times3$ unit matrix, and consider the
action-angle chart
$$
  \tilde\cC = (\tilde a,\tilde \a) := (Z^Ta,Z\a) : 
  V \to (Z^T\cA)\times\toro^n
$$
on $V$. Clearly, at the points of $J^{-1}(c)$, the representative
$\tilde\cF$ of $f$ in this chart system is independent of the
last six angles $\tilde\a$. In terms of the Fourier series
for $\pi:V\to \pi(V)=:B_V$ relative to the atlas formed by the
chart $\tilde\cC$ this means that all
harmonics $f^\nu$ with $\nu=(\mu,\eta)\in \interi^{n-6}
\times\interi^6$ and $\eta\not=0$ vanish at the points of $J^{-1}(c)$
and hence
\beq{Sp-locale}
  \Sp{f,b}\subseteq \interi^{n-6}\times \{0\} 
  \qquad
  \forall\ b\in\pi(J^{-1}(c))
  \,.
\eeq
We now extend this result to all of $B_V=\pi(V)$ (this is the key
point, and uses that $f$ has property FG2). \for{Sp-locale} implies
that, for all $\mu \in \interi^{n-6}$ and $\eta \in
\interi^{6}\setminus\{0\}$, $B_V^{(\mu,\eta)}(f) \cap
\pi(J^{-1}(c)) = \emptyset$ and hence
$$
   B_V\setminus B_V^{(\mu,\eta)}(f) \supset \pi(J^{-1}(c)) 
$$
(the sets $B_V^\nu(f)$ are defined as in \for{B-nu}, but with $B$
replaced by $B_V$). 
$\pi(J^{-1}(c))$ is diffeomorphic to $\{I\in\reali^{n-3}:(I,c)\in\cA\}$
and so is not countable. Since $f$ has property FG2 (so that each
$B_V\setminus B_V^{\nu}(f)$ either equals $B_V$ or is countable) this
implies $B_V^{(\mu,\eta)}(f) =\emptyset$ if $\eta\not=0$. Hence, the
representative of $f$ in the chart $\tilde\cC$ is
independent of the last 6 angles at all points of $V$. 
Note also that the actions of the chart $\tilde\cC$ take values in a
convex set (linear transformations map convex
sets to convex sets).

This shows that the original system, restricted to the subset
$V$, is invariant under a fully-Hamiltonian action of $\toro^6$, which
is given by translation of the last 6 angles of the chart $\tilde\cC$.
If all the reduced systems $(V_c,\s_c,f_c)$, $c\in\reali^6$, 
are symplectic, then the proof of Theorem \ref{thm2} is terminated,
with $k=6$. 

If not, this procedure can be repeated to obtain reduced systems with
$n-9$, $n-12$, $\dots$ degrees of freedom. The iteration stops when
all reduced systems are symplectic; this certainly happens when $0$,
$1$ or $2$ degrees of freedom are reached. 

\section{Examples}
\label{s:esempi}

\subsection{The reconstruction equation}
\label{ss:reconstruction1}
We provide now some simple examples of non-vertical fully-Hamiltonian
vector fields, and describe some of their properties. This requires
a reconstruction procedure from the symplectic reduced systems of
Theorem \ref{thm2} to the almost-symplectic unreduced one.

We assume $M=\renne\times\toro^n\ni(a,\a)$ and $\s=\saa$ as in
\for{eq:sigma}. Consider a function $f\in\cH_F(M)$ that has property
FG2. It follows from Theorem \ref{thm2} that $f$ is invariant under the
action of a $k$-dimensional subtorus $\toro^k$ of the torus $\toro^n$
that acts by translation of the angles $\a$, for some $3\le k\le n$.

When $k<n$, it follows from Remark i. at the end of Section \ref{ss:results}
that there is a matrix $Z\in SL_\pm(n,\interi)$ such that the representative of
$f$ in the new action-angle coordinates $(Za,Z^{-T}\a)$ (that we
still denote $f$) is independent of the last $k$ angles. These
new action-angle coordinates can be partitioned as
\beq{aa-10}
  (I,J,\phi,\psi) \in
  \reali^{r}\times\reali^{k} \times 
  \toro^{r}\times\toro^{k}
\eeq
with $r=n-k$ and
$$
   \der f\psi=0 \,.
$$
The action of $\toro^k$ acts by translations of the angles $\psi$ and
has momentum map $J$. 

We restrict to a level set
$J^{-1}(c)=\reali^{r}\times\toro^{r}\times \toro^k\ni
(I,\phi,\psi)$ of $J$. The restriction of the equations of motion
\for{queste} of $X_f^\s$ to such a level set is
\begin{equation}
\label{eq-1}
  \dot I = -\der {f} \phi(I,c,\phi) \,,\qquad 
  \dot\phi = \der {f} I(I,c,\phi)  + A_{II}(I,c) \der{f}\phi(I,c,\phi) 
  \,,\qquad (I,\phi)\in\reali^r\times\toro^r \,,
\end{equation}
and
\begin{equation}
\label{eq-2}
  \dot\psi = \der {f} J(I,c,\phi)  + A_{JI}(I,c)\der{f}\phi(I,c,\phi)
  \,,\qquad \psi\in\toro^k \,. 
\end{equation}

The coordinates $(I,\psi)$ serve as coordinates on the reduced space
$M_c=\reali^{r}\times\toro^{r}$. The symplectic structure $\s_c$ of a
reduced space $M_c$ takes the form
$$
  \s_c(I,\phi) = dI\wedge d\phi + \frac12 A_{II}(I,c) dI\wedge dI \,.
$$
Note that, unless $A_{II}=0$ (which is always the case if $r=1$,
namely, $k=n-1$), $(I,\phi)$ are not Darboux coordinates for $\s_c$. 
The reduced function $f_c(I,\phi)$ on $M_c$ equals the restriction
$f(I,c,\phi)$ of $f$ to $J^{-1}(c)$ and the reduced equations of
motion on $M_c$ are 
\begin{equation}
\label{eq-1}
  \dot I = -\der {f_c} \phi(I,\phi) \,,\qquad 
  \dot\phi = \der {f_c} I(I,\phi)  + A_{II}(I,c) \der{f_c}\phi(I,\phi) 
  \,,\qquad (I,\phi)\in\reali^r\times\toro^r\,,  
\end{equation}
and are symplectic-Hamiltonian (even though written in non-Darboux
coordinates).  

Equation \for{eq-2}, together with $J=c$, is the
reconstruction equation. The possibly non-symplectic nature of the
system is encoded in the terms $A_{JI}\der f{\phi}$.

\subsection{Examples}
\label{ss:esempi.2}

We describe now various possibilities met in the above reconstruction
process, and the ensuing dynamical properties of the fully-Hamiltonian systems.
This will also demonstrate the three statements after Theorem \ref{thm2}.
An important role is played by the dimension $r=n-k$ of the reduced
spaces.

{\it 1. $r=0$. } This case is met when the function $f$ is
independent of all angles and is basic. Thus, as seen in
Section \ref{ss:3.2}, $X_f^\s$ is vertical and Hamiltonian, and
completely integrable, with respect to the (in this example, globally
defined) symplectic structure $da\wedge d\a$.

We note that this argument could be generalized to prove to any
dimension $n\ge3$ the result of item ii. of Theorem \ref{thm1},
under however the stronger hypothesis that $f$ has
property FG2. 

{\it 2. $r=1$. } In the adapted system of action-angle coordinates
\for{aa-10}, $f$ depends on $(I,J)=(I,J_1,\ldots,J_{n-1})$ and on the
single angle $\phi\in\toro^1$. The argument used in the proof of Theorem
\ref{thm3} applies to this case as well (see the Remark at the end of
Section \ref{ss:ProofThm3}) and shows that $X_f^\s$ is Hamiltonian,
and completely integrable, with respect to a certain symplectic
structure $\tilde\s^\mathrm{aa}$ on $M$, see \for{sigma-3}. 

Specifically, the $n$ independent integrals in involution, relatively
to $\tilde\s^\mathrm{aa}$, are $f$ and the $n-1$ actions $J$ (unless
$\der f\phi=0$, in which case they are the $J$'s and $I$, as
discussed after \for{sigma-3}). If the common level sets of these
integrals are compact, motions $X_f^\s$ are quasi-periodic on tori of 
dimension $n$.

This generalizes to any dimension the result of Theorem
\ref{thm3} relative to the case $n=3$, under the stronger
hypothesis that $f$ has property FG2. 

One can say more if the regular level sets of a reduced Hamiltonian
$f_c$ are compact (hence, being one-dimensional, unions of closed
curves). The reduced dynamics in the reduced space $M_c$ consists of
periodic orbits, equilibria, and curves asymptotic to the equilibria.
Thus, with the exclusion of the latter, all reduced motions
reconstruct to quasi-periodic motions by known results on 
reconstruction from periodic dynamics
\cite{field1,field2,krupa,hermans}.

An example, with $n=4$, is the following. Consider three functions
$g_{2}, g_{3}, g_{4}:\reali^2\to\reali$ and the matrix
$$
  A_4=\left(
  \begin{matrix} 
     0 &-g_{2}(a_1,a_2) &-g_{3}(a_1,a_3)  &-g_{4}(a_1,a_4)  \cr
     g_{2}(a_1,a_2) &0 &a_4 &0 \cr
     g_{3}(a_1,a_3)   &-a_4 &0 &0 \cr
     g_{4}(a_1,a_4)   &0 &0 &0
  \end{matrix}
  \right) \,.
$$
The tensor $C$ associated to $A_4$, see \for{eq:C}, has all entries
whose indices are a permutation of $\{2,3,4\}$ equal to $\pm1$.
Therefore, for any choice of the functions $g_2,g_3,g_4$ the matrix
$A_4$ leads to a 2-form $\s$ on $M=\reali^ 4\times \toro^4$ which is not
symplectic. (The non-closedness of $\s$ comes only from the
lower-right $3\times 3$ diagonal block of $A_4$; the functions
$g_2,g_3,g_4$ affect only the reconstruction, hence the properties of
$X^\s_f$). 
For any function $f:\reali^4\times\toro^4\to\reali$, the quantities
$\sum_kC_{ijk}\frac{\partial f}{\alpha_k}$ are either $0$ or 
$\pm \der f{\a_2}$, $\pm \der f{\a_3}$, $\pm\der f{\a_4}$. 
Therefore, $X^\s_f$ is fully-Hamiltonian if and only if $f$ is
independent of the three angles $\a_2,\a_3,\a_4$. Symplectization of 
the system is obtained via reduction under the $\toro^3$-action given
by translations of the angles $(\a_2,\a_3,\a_4)$. In the notation of
Section \ref{ss:reconstruction1}, $I=a_1$, $\phi=\a_1$,
$J=(a_2,a_3,a_4)$ and $\psi=(\a_2,\a_3,\a_4)$. Thus, $X^\s_f$ is
not-vertical if
\[
   \der f {\a_1}\not=0 \,.
\]
For a value $c=(c_2,c_3,c_4)$ of $J$, the reduced Hamiltonian is
$f_c(a_1,\a_1)=f(a_1,c_2,c_3,c_4,\a_1)$ and the reduced equations are
$\dot a_1=-\der {f_c}{\a_1}(a_1,\a_1)$, $\dot \a_1=\der {f_c}
{a_1}(a_1,\a_1)$. Hence, the reduced vector field is Hamiltonian with
respect to the symplectic structure $dI\wedge d\phi$. 
The equations of motion of $X_f^\s$ are \for{queste}. They are
Hamilton equations relative to the symplectic form $da_1\wedge
d\a_1+\sum_{i=2}^4da_i\wedge d\chi_i$ with $\chi_i=\a_i-\int
g_i(a_1,a_i)da_1$ (see \for{quelle}). 

{\it 3. $r=2$. } In the previous two cases, all fully-Hamiltonian
vector fields are Hamiltonian, and completely integrable, with
respect to some symplectic structure. When $k=n-2$, the reduced
spaces have $2$ degrees of freedom and there is more freedom. (The
same happens if $k\le n-2$ as in point 4. below). 

On the one hand, in this case the fully-Hamiltonian systems need not
be Hamiltonian with respect to any symplectic structure, not even
semi-globally. (If there is more than one $I$, the argument used in
the proof of Theorem \ref{thm3} to build a symplectic structure
requires  some extra closedness conditions on the matrix $A_{IJ}$). 

On the other hand, from a dynamical point view, the reduced systems
are standard Hamiltonian systems with $2$ degrees of freedom and could be
chaotic rather than integrable. 

As an example with $n=5$, consider the matrix
$$
  A_5=\left(
  \begin{matrix} 
     0 &g_{12} &g_{13} &g_{14} &g_{15} \cr 
     -g_{12} &0 &g_{23} &g_{24} &g_{25} \cr
     -g_{13} &-g_{14} &0 &a_5 &0 \cr
     -g_{14} &-g_{24} &-a_5 &0 &0 \cr
     -g_{15} &-g_{25} &0 &0 &0 \cr
  \end{matrix}
  \right) 
$$
where each function $g_{ij}$ depends on the two actions $a_i$ and
$a_j$. Here too, any choice of these functions gives an
almost-symplectic form in $\reali^5\times\toro^5$ which is not
symplectic, and a partially-Hamiltonian vector field $X_f^\s$ is
fully-Hamiltonian if and only if $f$ is independent of
$\alpha_3,\alpha_4,\alpha_5$. It is non-vertical
if $\der f{\a_1}\not=0$ and $\der f{\a_2}\not=0$. 

An instance of a fully-Hamiltonian function, in this case, is the function
$$
  \frac{a_1^2}2+a_2-(1+\cos\alpha_2)\cos\alpha_1 +a_3+a_4+a_5 \,,
$$
which describes a system formed by a periodically perturbed pendulum,
which is well known to be nonintegrable and has chaotic dynamics, and
three uncoupled oscillators. Thus, there are non-vertical
fully-Hamiltonian vector fields with chaotic dynamics.

{\it 4. $r>2$. } Finally, we note that it could happen that the
symplectic reduced systems have more than 2 degrees of freedom. This
happens if,  in the iteration of the proof of Theorem \ref{thm2}, 
all almost-symplectic reduced spaces with $r>2$ degrees of freedom are
symplectic, so that the iteration stops. (In fact, the reduction
procedure {\it must} stop at that step because, in a symplectic
manifold, the condition of full-Hamiltonianity does not impose any
symmetry).

As a simple example, consider the almost-symplectic structure
$$
  \s = \sum_{i=1}^6da_i\wedge d\a_i
  + \frac 12 \sum_{i,j=1}^5 (A_5)_{ij}da_i\wedge da_j  
$$
on $\reali^6\times\toro^6$, where the matrix $A_5$ is as in item 3.
above. All functions independent of $\a_3,\a_4,\a_5$ are now
fully-Hamiltonian, but clearly, after reducing under translations of
these three angles, a symplectic system with three degrees of freedom
is met, and since the dependence on the angle $\a_6$ can be anything,
in general it cannot be further reduced under a $\toro^3$-action.\\

{\bf Acknowledgements.}
F.F. has been partially supported by the MIUR-PRIN project 20178CJA2B
{\it New Frontiers of Celestial Mechanics: theory and applications}.
Both authors have been partially supported by the MIUR-PRIN project
2022FPZEES {\it Stability in Hamiltonian dynamics and beyond}.
The authors are grateful to
Giovanna Carnovale, Andrea Giacobbe, Umberto Marconi and
James Montaldi for useful conversations on some technical points, and
to Heinz Han{\ss}mann for useful comments on a preliminary version.
The authors are members of {\it GNFM-INDAM}.

\end{document}